\documentclass[reqno, 10pt]{amsart}

\usepackage{color}
\usepackage{mathrsfs}
\usepackage{amssymb}

\newtheorem{thm}{Theorem}[section]
\newtheorem{cor}[thm]{Corollary}
\newtheorem{lem}[thm]{Lemma}
\newtheorem{prop}[thm]{Proposition}
\theoremstyle{definition}
\newtheorem{defn}[thm]{Definition}
\theoremstyle{remark}
\newtheorem{rem}[thm]{Remark}

\numberwithin{equation}{section}

\newcommand{\eqy}{\backsim}

\begin{document}

\title{Observability and nonlinear filtering}
\author{Ramon van Handel}
\address{Physical Measurement and Control 266-33 \\
California Institute of Technology \\
Pasadena, CA 91125 \\
USA}
\email{ramon@its.caltech.edu}
\thanks{This work was supported by the ARO}



\begin{abstract} 
This paper develops a connection between the asymptotic stability of
nonlinear filters and a notion of observability.  We consider a general
class of hidden Markov models in continuous time with compact signal
state space, and call such a model observable if no two initial measures
of the signal process give rise to the same law of the observation
process.  We demonstrate that observability implies stability of the
filter, i.e., the filtered estimates become insensitive to the initial
measure at large times.  For the special case where the signal is a
finite-state Markov process and the observations are of the white noise
type, a complete (necessary and sufficient) characterization of filter
stability is obtained in terms of a slightly weaker detectability
condition.  In addition to observability, the role of controllability in 
filter stability is explored.  Finally, the results are partially extended 
to non-compact signal state spaces. 
\end{abstract}

\maketitle

\section{Introduction}

Consider the deterministic linear control system
\begin{equation*}
\begin{split}
	\frac{d}{dt}\,x(t) &= A\,x(t) + \Sigma\,u(t), \\
	y(t) &= C\,x(t),
\end{split}
\end{equation*}
where $x(t)$ is the system state, $u(t)$ is the control input and $y(t)$
is the observation signal.  Such a system is called \textit{observable}
if there exist no $x\ne x'$ such that $\{y(t):t\ge 0\}$ is the same when
$x(0)=x$ and $x(0)=x'$ (for any control $u(t)$), and is called
\textit{controllable} if for any $x,x'$ and $t>0$, there is a control signal
$u(t)$ such that the solution with $x(0)=x$ satisfies $x(t)=x'$.  It is
well known \cite{BJ68,OP96} that observability and controllability are 
intimately related with the asymptotic properties of the conditional estimates
in the linear filtering problem
\begin{equation*}
\begin{split}
	dX_t &= AX_t\,dt + \Sigma\,dW_t, \\
	dY_t &= CX_t\,dt + dB_t,
\end{split}
\end{equation*}
where $W_t$ and $B_t$ are independent standard Wiener processes.  In 
particular, it is found that the filtered estimates become insensitive
to the law of $X_0$ at large times, i.e.,\footnote{
	Here $\mathbf{P}^\mu$ is the law under which $X_0\sim\mu$,
	$\mathscr{F}_t^Y=\sigma\{Y_s:s\le t\}$, and the versions of the
	conditional expectations are chosen to coincide with those
	computed by the Kalman filter.}
$|\mathbf{E}^\mu(f(X_t)|\mathscr{F}^Y_t)-\mathbf{E}^\nu(f(X_t)|
\mathscr{F}_t^Y)|\to 0$ as $t\to\infty$ for any pair of initial laws 
$X_0\sim\mu,\nu$, whenever the associated linear control system is observable
and controllable \cite{OP96}.  This is called the \textit{stability} property
of the filtering problem, and is of significant importance from the practical
point of view as it ensures the robustness of the filtered estimates with 
respect to modelling errors and approximations.  The purpose of this paper
is to demonstrate that the connection between observability, contollability,
and the linear filtering problem has a natural counterpart in a large variety
of nonlinear filtering problems for Markovian signal-observation models.

We consider a signal process $X_t$ and observation process $Y_t$ in
continuous time, where both $X_t$ and $(X_t,Y_t)$ are assumed to be
Markov processes and $Y_t$ is assumed to satisfy a mild condition which
ensures, in essence, that the observation noise is memoryless (the most
common observation models in continuous time, additive white noise
observations and counting observations, both satisfy this requirement). 
In addition, we will mostly assume that the signal process takes values
in a compact state space.  The significance of the compactness
assumption and some extensions of the results to the non-compact case
are discussed in section \ref{sec:non-compact}.

In this general setting, the model is called \textit{observable} if
there do not exist initial measures $\mu\ne\nu$ such that $\{Y_t:t\ge
0\}$ has the same law when $X_0\sim\mu$ and $X_0\sim\nu$.  One of the
main results of this paper (corollary \ref{cor:stab}) is that if the
model is observable, then
$|\mathbf{E}^\mu(f(X_t)|\mathscr{F}^Y_t)-\mathbf{E}^\nu(f(X_t)|
\mathscr{F}_t^Y)|\to 0$ as $t\to\infty$ $\mathbf{P}^\mu$-a.s.\ for any
bounded continuous $f$ and any pair of absolutely continuous initial
laws $\mu\ll\nu$.

On the other hand, the notion of \textit{controllability} is replaced by
a certain regularity property of the signal transition probabilities
(see definition \ref{def:regular}). Under the additional assumption that
the observations are of the white noise type and nondegenerate, we show
that observability and regularity of the signal imply that
$|\mathbf{E}^\mu(f(X_t)|\mathscr{F}^Y_t)-\mathbf{E}^\nu(f(X_t)|
\mathscr{F}_t^Y)|\to 0$ as $t\to\infty$ $\mathbf{P}^\mu$-a.s.\ for any
bounded continuous $f$ and initial laws $\mu,\nu$ (corollary
\ref{cor:genultra}), without the absolute continuity requirement.  When
the signal is the solution of a stochastic differential equation, then
the regularity property is directly related to the controllability of an
associated deterministic control system.  This result is thus entirely
parallel to the observability-controllability criterion for the
stability of the Kalman filter.

A natural test case for the general theory is the setting where the
signal state space is a finite set.  The simplicity of this setting
allows a particularly transparent insight into the nature of the
observability and controllability properties, while this setting is also
of particular practical importance due to the fact that the associated
nonlinear filters are finite dimensionally computable.  By combining the
general results in this paper with known filter stability results for
ergodic signals \cite{BCL}, we find that a complete characterization of
filter stability is possible for finite state signals with nondegenerate
white noise type observations.  In particular, we will find necessary
and sufficient criteria for stability (theorems \ref{thm:detfin} and
\ref{thm:ultdetfin}) which can be verified explicitly for a given model
through straightforward linear algebra techniques.  The fact that such a
complete characterization is possible, albeit in a particularly simple
case, suggests that the notion of observability used in this paper is in
some sense a fundamental ingredient of the filter stability problem.

The stability of nonlinear filters has been studied actively in the last
few years following the pioneering contributions of Ocone and Pardoux
\cite{OP96} and Zeitouni {\it et al.}\ \cite{delyon,AZ97}.  An excellent
overview of previous work and an extensive list of references can be
found in \cite{ChSurv}.  The majority of the results on this topic
assume that the signal is an ergodic Markov process and that the
observation process is of the white noise type.  Such results are
complementary to the results obtained in this paper; indeed, for the
complete characterization of stability in the finite state setting, it
is essential to combine our results with results that are specific to
ergodic signals.  Moreover, ergodic results may fail to hold true for
certain degenerate observation models (see the counterexample in
\cite{BCL}), while the results of this paper still hold in that setting.
On the other hand, much is known about the rate of convergence of
differently initialized filters in particular cases (see, e.g.,
\cite{delyon,AZ97,BCL}), while our results do not provide such
information.  In a much more general setting, some results on the
asymptotic properties of nonlinear filtering errors can be found in Clark, 
Ocone and Coumarbatch \cite{COC99}.  Their point of view is close to the one
used in this paper, but their results do not establish stability of the
filter.

The approach in this paper is inspired by the observation of Chigansky
and Liptser \cite{CL} that, by virtue of a martingale convergence
argument, certain predictive estimates of the observations are always
stable.  This implies, in particular, that the filtered estimates of the
signal are stable for a particular class of functions, which can be
written as the conditional mean of a functional of the observation
process given the initial value of the signal (see section
\ref{sec:filstob}).  The heart of the argument that leads to the
observability criterion is the characterization of this class of
functions (proposition \ref{prop:dense}), which is achieved using a
corollary of the Hahn-Banach theorem.

The remainder of this paper is set up as follows.  In section
\ref{sec:model}, we introduce the general signal-observation model that
will be used in most of the paper, and fix the notation for the
remainder of the paper.  Section \ref{sec:spaces} is devoted to the
study of the notion of observability and its connection to the class of
functions that can be obtained as predictive estimates of the
observations.  Section \ref{sec:filstob} connects these concepts to the
stability of the filter. In section \ref{sec:obsmodels}, we show how the
notion of observability can be characterized in the common cases of
white noise type and counting observations, and we find a particularly
simple sufficient criterion for observability in those cases.  Section
\ref{sec:fss} is devoted to the finite state case; an explicit criterion
is found for observability, and stability is completely characterized
for absolutely continuous initial measures. Section \ref{sec:ultra}
explores the connection between controllability, regularity, and the
stability of filters for arbitrary initial conditions; a complete
characterization is again given for the finite state case. Finally,
section \ref{sec:non-compact} discusses the significance of the
compactness assumption made in the previous sections, and provides some
partial extensions of previous results.  A simple but apparently
unknown result for the Kalman filter is briefly discussed in the
appendix.

Before proceding to the main part of the paper, we make a few general remarks.

\begin{rem}
The main results of this paper can be adapted in a straightforward fashion
to the discrete time setting.
\end{rem}

\begin{rem}
It is not difficult to show that when our general notions of
observability and controllability are applied to the linear filtering
model, one obtains precisely the classical observability-controllability
criteria for the Kalman filter.  Unfortunately, the results in section
\ref{sec:non-compact} for the non-compact case are not sufficiently
powerful to recover the stability of the Kalman filter, except then the
signal itself is asymptotically stable (i.e., the matrix $A$ has only
eigenvalues with strictly negative real parts).  The latter case is not
particularly interesting for the Kalman filter, as the more general
detectability criterion (see \cite{OP96} and the appendix) makes
observability irrelevant in this setting.  The fact that the Kalman filter 
with an unstable signal is not covered is a major shortcoming of the
results in this paper.
\end{rem}

\begin{rem}
With the exception of the Kalman filter, the finite state case, and
observations with an invertible observation function (lemma
\ref{lem:onetoone}), the observability property appears to be 
difficult to verify for a given model.  For practical applications, it
is thus necessary to develop explicitly verifiable sufficient criteria for
observability (see section \ref{sec:ultrageneral} for further
discussion).  
\end{rem}

\section{The signal-observation model}
\label{sec:model}

The goal of this section is to set up the model for the signal and
observation processes, and to fix the notation that will be used in the
following.

Let us begin by introducing the basic objects that make up the model.
\begin{enumerate}
\item The \textit{signal state space} $\mathbb{S}$ is a compact Polish space.
\item The \textit{observation state space} $\mathbb{O}=\mathbb{R}^p$ for
some $p<\infty$.
\item The \textit{signal-observation process} $(X_t,Y_t)_{t\in[0,\infty[}$ 
is a time-homogeneous $\mathbb{S}\times\mathbb{O}$-valued Feller-Markov process
with c{\`a}dl{\`a}g paths.
\item The \textit{signal process} $(X_t)_{t\in[0,\infty[}$ is a Feller-Markov
process in its own right.
\item The \textit{observation process} $(Y_t)_{t\in[0,\infty[}$ has
conditionally independent increments given the signal process
$(X_t)_{t\in[0,\infty[}$, and $Y_0=0$.
\end{enumerate}
This can be viewed as a hidden Markov model in continuous time,
where $Y_t$ is the observable component and $X_t$ is the nonobservable 
component.

For any locally compact Polish space $S$, we denote by $\mathscr{B}(S)$
the Borel $\sigma$-algebra, by $\mathcal{C}(S)$ the space of continuous
functions, by $\mathcal{C}_b(S)$ the space of bounded continuous
functions, by $\mathcal{C}_0(S)$ the space of continuous functions that
vanish at infinity, by $\mathcal{M}(S)$ the space of finite signed
measures on $\mathscr{B}(S)$, by $\mathcal{P}(S)$ the space of
probability measures on $\mathscr{B}(S)$, and by $\mathcal{M}_c(S)$
($\mathcal{P}_c(S)$) the finite signed (probability) measures with
compact support. Note that when $S$ is compact,
$\mathcal{C}(S)=\mathcal{C}_b(S)= \mathcal{C}_0(S)$.

It is convenient to construct the signal-observation process
$(X_t,Y_t)_{[0,\infty[}$ on its canonical probability space.  To this
end, define $\Omega^X=D([0,\infty[\mbox{}; \mathbb{S})$ and
$\Omega^Y=D([0,\infty[\mbox{};\mathbb{O})$, i.e., $\Omega^X$ and
$\Omega^Y$ are the spaces of $\mathbb{S}$-valued and $\mathbb{O}$-valued
c{\`a}dl{\`a}g paths, endowed with the Skorokhod topology. We will work
on the probability space $\Omega=\Omega^X\times\Omega^Y$, equipped with
its Borel $\sigma$-algebra
$\mathscr{F}=\mathscr{B}(\Omega^X\times\Omega^Y)$, and choose
$X_t:\Omega\to\mathbb{S}$ and $Y_t:\Omega\to\mathbb{O}$ to be the
canonical processes $X_t(x,y)=x(t)$ and $Y_t(x,y)=y(t)$.  Furthermore,
we define the natural filtrations 
$$
        \mathscr{F}_t^{X}=\sigma\{X_s:s\le t\},\qquad
        \mathscr{F}_t^{Y}=\sigma\{Y_s:s\le t\},\qquad
	\mathscr{F}_t = \sigma\{(X_s,Y_s):s\le t\},
$$
and the filtration generated by the observation increments
$$
	\mathscr{G}_t^Y=\sigma\{Y_s-Y_0:s\le t\}.
$$
We will denote $\mathscr{F}^X=\mathscr{F}^X_\infty=\bigvee_{t\ge 0}
\mathscr{F}_t^X$, and we define $\mathscr{F}^Y$ and $\mathscr{G}^Y$
similarly.

Let $T_t:\mathcal{C}_0(\mathbb{S}\times\mathbb{O})\to
\mathcal{C}_0(\mathbb{S}\times\mathbb{O})$ and $P_t(x,y,A)$
($t\in[0,\infty[\mbox{}$, $A\in\mathscr{B}(\mathbb{S}\times\mathbb{O})$)
be the Markov semigroup of the signal-observation process and the
associated transition probabilities.  By the Feller assumption, we can
construct a process with c{\`a}dl{\`a}g paths which possesses the
desired transition probabilities \cite[theorem 17.15]{Kal97}.  Hence there 
exists a family of probability measures
$\{\mathbf{P}_{(x,y)}:(x,y)\in\mathbb{S}\times\mathbb{O}\}$ on
$(\Omega,\mathscr{F})$ such that for every $(x,y)$, the process
$(X_t,Y_t)$ is a Markov process with respect to the filtration
$\mathscr{F}_{t}$ under $\mathbf{P}_{(x,y)}$ with transition
probabilities $P_t(x,y,A)$ and initial law
$(X_0,Y_0)\sim\delta_{\{(x,y)\}}$, and
$(x,y)\mapsto\mathbf{P}_{(x,y)}(A)$ is measurable for every
$A\in\mathscr{F}$.  In particular, under the probability measure
$$
        \mathbf{P}_\mu(A) = \int_{\mathbb{S}\times\mathbb{O}}
                \mathbf{P}_{(x,y)}(A)\,\mu(dx,dy),\qquad
        A\in\mathscr{F},\quad
        \mu\in\mathcal{P}(\mathbb{S}\times\mathbb{O}),
$$
$(X_t,Y_t)$ is a Markov process with respect to the filtration 
$\mathscr{F}_t$ with transition probabilities $P_t(x,y,A)$ and initial law
$(X_0,Y_0)\sim\mu$.  We recall that the Markov property can be expressed as 
follows \cite[proposition III.1.7]{RY98}: for bounded
$\mathscr{F}$-measurable $\xi$
$$
        \mathbf{E}_\mu(\xi\circ\theta_t|\mathscr{F}_t) =
        \mathbf{E}_{(X_t,Y_t)}(\xi)\quad\mathbf{P}_\mu\mbox{-a.s.}
	\qquad
        \mbox{for all }t>0,
$$
where $\mathbf{E}_\mu$, $\mathbf{E}_{(x,y)}$ denote the expectations
with respect to the measures $\mathbf{P}_\mu$ and $\mathbf{P}_{(x,y)}$,
and $\theta_t:\Omega\to \Omega$ is the canonical shift
$\theta_t(x,y)(s)=(x(s+t),y(s+t))$.

It is convenient, without loss of generality, to replace the various
$\sigma$-algebras and filtrations defined above by their usual
augmentations with respect to the family $\{\mathbf{P}_\mu:\mu\in
\mathcal{P}(\mathbb{S}\times\mathbb{O})\}$ \cite[section 1.4]{RY98}, and
we will make this replacement from this point onwards.  A significant
advantage of this choice is that if a bounded process $Z_t$ has
c{\`a}dl{\`a}g paths, and the filtration $\mathscr{G}_t$ satisfies the
usual conditions, then we can choose a version of
$\mathbf{E}(Z_t|\mathscr{G}_t)$, for every time $t$, so that the process
$t\mapsto\mathbf{E}(Z_t|\mathscr{G}_t)$ has c{\`a}dl{\`a}g paths
\cite[chapter VI, theorem 47]{DM}, \cite[theorem 6]{rao}. \textit{In the
following, whenever such processes are encountered, their c{\`a}dl{\`a}g
versions are always implied.}

Finally, let us make precise the conditions on the signal and
observations, i.e., that the signal is a Markov process in its own right
and that the observation process has conditionally independent
increments given the signal process.  Both these properties can be
simultaneously introduced through the following requirement.
\begin{itemize}
\item The signal is a Markov process in its own right, and
	the observation process has conditionally independent
	increments given the signal process, in the following sense:
	if the random variable $\xi$ is bounded and $\mathscr{F}^X\vee
	\mathscr{G}^Y$-measurable, then the map
	$(x,y)\mapsto \mathbf{E}_{(x,y)}(\xi)$ does not depend on $y$.
\end{itemize}
Using the Markov property of $(X_t,Y_t)$, this implies that
$\mathbf{E}_\mu(\xi|\mathscr{F}_s)=\mathbf{E}_\mu(\xi|X_s)$
$\mathbf{P}_\mu$-a.s.\ whenever $\xi$ is $\sigma\{X_t:t\ge
s\}$-measurable, which establishes that $X_t$ is an
$\mathscr{F}_t$-Markov process as desired.  On the other hand, we find
that for any bounded, $\sigma\{Y_{t+s}-Y_s:t>0\}$-measurable random
variable $\xi$, there exists a measurable function
$f:\mathbb{S}\to\mathbb{R}$ such that
$\mathbf{E}_\mu(\xi|\mathscr{F}_s)=f(X_s)$ $\mathbf{P}_\mu$-a.s.\ for
any initial measure $\mu$ (by the Markov property).  This expresses the
fact that the additional randomness introduced by the observation
process is memoryless. As we will see, the two most common types of
observations encountered in continuous time problems, white noise type 
observations and counting observations, satisfy this property.

It remains to note that the assumption $Y_0=0$ means that we will be
interested in initial measures of the form $\mu\times\delta_{\{0\}}\in
\mathcal{P}(\mathbb{S}\times\mathbb{O})$, where
$\mu\in\mathcal{P}(\mathbb{S})$. We therefore introduce the following
notation: for any $\mu\in\mathcal{P}(\mathbb{S})$, we define
$\mathbf{P}^\mu=\mathbf{P}_{\mu\times\delta_{\{0\}}}$. Similarly,
$\mathbf{E}^\mu$ denotes the expectation with respect to
$\mathbf{P}^\mu$.

\begin{rem}
There is no loss of generality in assuming that $Y_0=0$.  Indeed, consider
an arbitrary initial measure $\mu\in\mathcal{P}(\mathbb{S}\times\mathbb{O})$.
Then by \cite[lemma 2.4]{COC99}
$$
	\mathbf{E}_\mu(f(X_t)|\mathscr{F}_t^Y) =
	\mathbf{E}_{\mu(\,\cdot\,|Y_0)\times\delta_{\{Y_0\}}}
	(f(X_t)|\mathscr{F}_t^Y),
$$
where $\mu(\,\cdot\,|Y_0)$ is a regular conditional probability of $X_0$ with 
respect to $Y_0$ under $\mu$.  But note that under any initial measure
of the form $\nu\times\delta_{\{a\}}$, our assumptions imply that
$\mathscr{F}^X\vee\mathscr{G}^Y$ is independent of
$\mathscr{F}_0^Y$, so that
$$
	\mathbf{E}_\mu(f(X_t)|\mathscr{F}_t^Y) =
	\mathbf{E}_{\mu(\,\cdot\,|Y_0)\times\delta_{\{Y_0\}}}
	(f(X_t)|\mathscr{G}_t^Y) = 
	\mathbf{E}^{\mu(\,\cdot\,|Y_0)}
	(f(X_t)|\mathscr{F}_t^Y)
$$
provided that we choose an appropriate version of the latter conditional
expectation that is defined $\mathbf{P}_\mu$-a.s.  Thus it suffices to 
consider the case $Y_0=0$.
\end{rem}

\section{Spaces of observable functions and nonobservable measures}
\label{sec:spaces}

Broadly speaking, the goal of this section is to investigate the following
question: what is the relation between the law of $X_0$ and the law of
$\mathscr{F}^Y$?  In the next section, we will see that this question has
immediate consequences for filter stability.

\begin{defn}
For $\mu,\nu\in\mathcal{P}(\mathbb{S})$, we write $\mu\eqy\nu$
whenever $\mathbf{P}^\mu|_{\mathscr{F}^Y}=\mathbf{P}^\nu|_{\mathscr{F}^Y}$.
In particular, $\backsim$ defines an equivalence relation on 
$\mathcal{P}(\mathbb{S})$.
\end{defn}

In words, if $\mu\eqy\nu$, then whenever $X_0$ has the law $\mu$ or $\nu$,
the same law of the observation process if obtained.  In particular,
no amount of statistics gathered from the observation process will allow
us to distinguish between $X_0\sim\mu$ and $X_0\sim\nu$.  This motivates
the following notion of observability, which is reminiscent (at least
in spirit) of the notion of observability used in linear systems
theory.

\begin{defn}
The filtering model is called \textit{observable} if $\mu\eqy\nu$
implies $\mu=\nu$.
\end{defn}

The following definition is key (we use the notation
$\mu(f)=\int f(x)\,\mu(dx)$).

\begin{defn}
Define the \textit{space of nonobservable measures} $\mathcal{N}$ as
$$
	\mathcal{N}=\{\alpha\mu_1-\alpha\mu_2\in\mathcal{M}(\mathbb{S}):
	\alpha\in\mathbb{R},~\mu_1,\mu_2\in\mathcal{P}(\mathbb{S}),~
	\mu_1\eqy\mu_2\}.
$$
Moreover, we define the \textit{space of observable functions}
$\mathcal{O}$ as
$$
	\mathcal{O}=\{f\in\mathcal{C}_b(\mathbb{S}):
		\mu_1(f)=\mu_2(f)\mbox{ for all }\mu_1\eqy\mu_2\}.
$$
\end{defn}
We attach to the nonobservable space $\mathcal{N}$ the following
intuitive interpretation: if we perturb the initial measure $X_0\sim\mu$
in the direction $\delta\in\mathcal{N}$ ($\mu\mapsto\mu+\delta$,
provided $\mu+\delta$ is again a probability measure), then the law of
the observation process does not change. The observable space
$\mathcal{O}$ then consists of those functions $f$ such that the
expectation of $f(X_0)$ is completely determined by the law of the
observation process.  Note that the filtering model is observable if and
only if every continuous function is observable, i.e., $\mathcal{O}=
\mathcal{C}_b(\mathbb{S})$, or, equivalently, if no nontrivial signed
measure is nonobservable, i.e., $\mathcal{N}=\{0\}$.

Our goal is to characterize the space $\mathcal{O}$.
Before we do this, let us recall a simple functional analytic device
which will be needed below \cite[chapter 4]{Rud73}.  Let $B$ be a Banach
space and denote by $B^*$ its topological dual.  Consider two (not
necessarily closed) linear subspaces $M\subset B$ and $N\subset B^*$. 
\begin{defn}
The \textit{annihilator} $M^\perp\subset B^*$ of $M$ is defined as
$$
        M^\perp=\{x^*\in B^*:\langle x^*,x\rangle=0\mbox{ for all }
                x\in M\}.
$$
Similarly, the annihilator $N^\perp\subset B$ of $N$ is defined as
$$
        N^\perp=\{x\in B:\langle x^*,x\rangle=0\mbox{ for all }
                x^*\in N\}.
$$
\end{defn}
The proof of the following lemma \cite[theorem 4.7]{Rud73} follows from a
straightforward application of the Hahn-Banach theorem.
\begin{lem}
\label{lem:annihilator}
$(M^\perp)^\perp=\overline{M}$, where $\overline{M}$ is the (norm-)closure of 
$M$ in $B$.
\end{lem}

Recall that $\mathcal{M}(\mathbb{S})$ is the topological dual of
$\mathcal{C}_b(\mathbb{S})$ by the Riesz-Markov theorem.  It is thus
easily verified from the definitions that
$\mathcal{O}=\mathcal{N}^\perp$.  What we will show is that there is a
dense subset $\mathcal{O}^0\subset\mathcal{O}$ such that every
$f\in\mathcal{O}^0$ can be written as $f(x) = \mathbf{E}_{(x,y)}(\xi)$
for some bounded $\mathscr{G}^Y$-measurable random variable $\xi$.

\begin{prop}
\label{prop:dense}
Let $\mathcal{O}^0$ be the linear span of functions of the form
$$
        \mathbf{E}_{(x,y)}(f_1(Y_{t_1}-Y_0)f_2(Y_{t_2}-Y_0)\cdots
        f_n(Y_{t_n}-Y_0)),
$$
for all $n<\infty$, $t_i\in D$ and bounded
continuous functions $f_i$ on $\mathbb{O}$, where $D$ is a dense subset of 
$[0,\infty[\mbox{}$.  Then $\mathcal{O}^0$ is dense in $\mathcal{O}$.  In 
particular, for any observable function $f\in\mathcal{O}$, there is a sequence
of functions $f_n\in\mathcal{O}^0$ such that $\|f-f_n\|\to 0$.
\end{prop}

\begin{proof}
By our assumptions, any $f\in\mathcal{O}^0$ only depends on
$\mathbb{S}$, and we find
\begin{multline*}
        \mathbf{E}_{(x,y)}(f_1(Y_{t_1}-Y_0)f_2(Y_{t_2}-Y_0)\cdots
        f_n(Y_{t_n}-Y_0)) = \\
        \mathbf{E}_{(x,0)}(f_1(Y_{t_1}-Y_0)f_2(Y_{t_2}-Y_0)\cdots
        f_n(Y_{t_n}-Y_0)) = \\
        \mathbf{E}_{(x,0)}(f_1(Y_{t_1})f_2(Y_{t_2})\cdots
        f_n(Y_{t_n})).
\end{multline*}
We claim that $\mathcal{O}^0\subset\mathcal{C}_b(\mathbb{S})$.  As
$\mathbb{S}$ is compact, it suffices to show that any
$f\in\mathcal{O}^0$ is continuous.  But if $x_n\to x$ in $\mathbb{S}$,
then by \cite[theorem 17.25]{Kal97} the measures $\mathbf{P}_{(x_n,0)}$
converge weakly to $\mathbf{P}_{(x,0)}$, and this in turn implies weak
convergence of the finite dimensional distributions on some dense subset
of times $D$ \cite[theorem 3.7.8]{EK86}.  The continuity of
$f\in\mathcal{O}^0$ follows directly from the previous expression.

To show that $\mathcal{O}^0$ is dense in $\mathcal{O}$, it suffices to
show that $(\mathcal{O}^0)^\perp=\mathcal{N}$ by lemma
\ref{lem:annihilator}.  Note that by \cite[theorem 16.6]{Bil} the finite
dimensional distributions in a dense set of times form a separating
class for probability measures on $D([0,\infty[\mbox{};\mathbb{O})$. 
Hence a standard monotone class argument shows that $\mu_1\eqy\mu_2$ if
and only if $\mathbf{E}^{\mu_1}(f_1(Y_{t_1})f_2(Y_{t_2})\cdots
f_n(Y_{t_n})) = \mathbf{E}^{\mu_2}(f_1(Y_{t_1})f_2(Y_{t_2})\cdots
f_n(Y_{t_n}))$ for all finite sets of times $t_i\in D$ and bounded
continuous $f_i$.  But using the previous equation display, this is
clearly the case if and only if $\mu_1(f)=\mu_2(f)$ for all
$f\in\mathcal{O}^0$.  Hence we find immediately that
$\mathcal{N}\subset(\mathcal{O}^0)^\perp$.  On the other hand, choose
any $\mu\in(\mathcal{O}^0)^\perp$ (with $\mu\ne 0$), and define
$\mu_1=\mu^+/\alpha$, $\mu_2=\mu^-/\alpha$ with
$\alpha=\mu^+(\mathbb{S})$ (here $\mu=\mu^+-\mu^-$ is the Hahn
decomposition of $\mu$).  Note that $\mu_1$ and $\mu_2$ are both
probability measures (due to the fact that $1\in\mathcal{O}^0$ implies
$\mu(\mathbb{S})=0$), and $\mu=\alpha\mu_1-\alpha\mu_2$.  But
$\mu_1(f)=\mu_2(f)$ for all $f\in\mathcal{O}^0$ implies
$\mu_1\eqy\mu_2$, so evidently $\mu\in\mathcal{N}$. Hence we have
established the converse inclusion
$\mathcal{N}\supset(\mathcal{O}^0)^\perp$, and the proof is complete. 
\end{proof}

\begin{rem}
One might hope that \textit{any} observable function $f\in\mathcal{O}$
can be written as $f(x)=\mathbf{E}_{(x,y)}(\xi)$ for some bounded
$\mathscr{G}^Y$-measurable $\xi$.  This seemingly plausible conjecture
need not hold true, however, as the following simplified example illustrates. 
Let $X$ be a $[0,1]$-valued random variable with law $\mu$, and let
$Y=X+\xi$ where $\xi$ is Gaussian with zero mean and unit variance.
Denote by $\mathbf{P}_\mu$ the joint law of $X$ and $Y$.  Then the same
argument used in the previous proof shows that any continuous
function $f:[0,1]\to\mathbb{R}$ can be written as the uniform limit of 
functions of the form $f_n(x) = \mathbf{E}_{\delta_{\{x\}}}(g_n(Y))$.  
However, any $f_n$ will necessarily be a smooth function (being the 
convolution of the bounded function $g_n$ with the Gaussian density), so 
that evidently not all $f$ can be expressed in this form.  Thus in general, 
an approximation result is the best one could hope for.
\end{rem}

\section{Filter stability and observability}
\label{sec:filstob}

We now connect the notions of observability introduced in
the previous section to the stability of the nonlinear filter.  Recall
that we are interested in determining, given a pair of initial measures
$\mu,\nu$, whether $\mathbf{E}^\mu(f(X_t)|\mathscr{F}_t^Y)$ and
$\mathbf{E}^\nu(f(X_t)|\mathscr{F}_t^Y)$ are close to each other for
large times $t$.  We will see that this is always the case when the
function $f$ is observable, i.e., when $f\in\mathcal{O}$, provided that
$\mu\ll\nu$.

The following lemma, which is inspired by a result of Chigansky and Liptser
\cite[theorem 2.1]{CL}, contains the essence of the convergence
argument.

\begin{lem}
\label{lem:cvgsimple}
Let $\mu,\nu\in\mathcal{P}(\mathbb{S})$ satisfy $\mu\ll\nu$.
Moreover, let $(\xi_t)_{t\in[0,\infty[}$ be $\mathscr{F}^Y$-measurable
random variables with $|\xi_t|\le K<\infty$ for all $t$, such that
the sample paths $t\mapsto\xi_t$ are c{\`a}dl{\`a}g.  Then we have
$$
        |\mathbf{E}^\mu(\xi_t|\mathscr{F}_t^Y)-
        \mathbf{E}^\nu(\xi_t|\mathscr{F}_t^Y)|
        \xrightarrow{t\to\infty}0\quad\mathbf{P}^\mu\mbox{-a.s.}
$$
\end{lem}

\begin{rem}
Recall that whenever conditional expectations are encountered, the
corresponding c{\`a}dl{\`a}g versions are implied.  Throughout the
following proofs, we will use the usual properties of conditional
expectations to obtain equalities and inequalities that, for every time
$t$, hold for all $\omega\in\Omega\backslash N_t$ where $N_t$ is a
$\mathbf{P}^\mu$-null set.  Because all the processes are
c{\`a}dl{\`a}g, however, the null set can be chosen independent of time
$t$, so that these equalities and inequalities hold for all $t$
simultaneously with unit probability.  We will use this fact below
without further comment.
\end{rem}

\begin{proof}
We begin by noting that \cite[lemma 2.1]{COC99}
$$
        \frac{d\mathbf{P}^\mu}{d\mathbf{P}^\nu} =
        \frac{d\mu}{d\nu}(X_0).
$$
By the Bayes formula, we obtain $\mathbf{P}^\mu$-a.s.
$$
        \mathbf{E}^\nu\left(\left.
                \frac{d\mu}{d\nu}(X_0)\,\right|
                \mathscr{F}_t^Y\right)
        \mathbf{E}^\mu(
                \xi_t|
                \mathscr{F}_t^Y) =
        \mathbf{E}^\nu\left(\left.
                \frac{d\mu}{d\nu}(X_0)~
                \xi_t\,\right|
                \mathscr{F}_t^Y\right).
$$
Introduce the notation
$$
        \varrho_t = \mathbf{E}^\nu\left(\left.
                \frac{d\mu}{d\nu}(X_0)\,\right|
                \mathscr{F}_t^Y\right),\qquad
        \varrho_\infty = \mathbf{E}^\nu\left(\left.
                \frac{d\mu}{d\nu}(X_0)\,\right|
                \mathscr{F}^Y\right).
$$
Then we find, using the fact that $\xi_t$ is $\mathscr{F}^Y$-measurable,
$$
        \varrho_t\,|\mathbf{E}^\mu(\xi_t|\mathscr{F}_t^Y)-
                \mathbf{E}^\nu(\xi_t|\mathscr{F}_t^Y)| =
        |\mathbf{E}^\nu((\varrho_\infty-\varrho_t)\,
                \xi_t|\mathscr{F}_t^Y)|\le
        K\,\mathbf{E}^\nu(|\varrho_\infty-\varrho_t|~
                |\mathscr{F}_t^Y).
$$
That this expression converges to zero $\mathbf{P}^\mu$-a.s.\ is established
in lemma \ref{lem:hunt} below.  But as $\varrho_t\to\varrho_\infty$
$\mathbf{P}^\nu$-a.s.\ by L{\'e}vy's upward theorem, we conclude the
convergence
$|\mathbf{E}^\mu(\xi_t|\mathscr{F}_t^Y)-\mathbf{E}^\nu(\xi_t|\mathscr{F}_t^Y)|
\to 0$ as $t\to\infty$ on $\{\omega\in\Omega:\varrho_\infty(\omega)>0\}\in
\mathscr{F}^Y$ (modulo a $\mathbf{P}^\mu$-null set), and the latter
set has $\mathbf{P}^\mu$-measure one.
\end{proof}

The proof of the previous lemma is not yet complete, as we still need to
show that
$\mathbf{E}^\nu(|\varrho_\infty-\varrho_t|~|\mathscr{F}_t^Y)\to 0$.  If
we were interested in $L^1$ convergence rather than a.s.\ convergence,
the result is trivially established.  Proving a.s.\ convergence would
appear to be a matter of applying Hunt's lemma \cite[chapter V, theorem
45]{DM}, whose proof is easily adapted to the continuous time setting. 
Unfortunately, this would require $\varrho_t$ to be dominated by an
integrable random variable, which may not be the case (to guarantee that
this is the case we could impose, e.g., a finite relative entropy
condition $D(\mu||\nu)<\infty$, see \cite[chapter V, sec.\ 25(c)]{DM}).  
Instead, we proceed by adapting Rao's proof of Hunt's lemma \cite[lemma 
2]{rao} to our setting.

\begin{lem}
\label{lem:hunt}
$\mathbf{E}^\nu(|\varrho_\infty-\varrho_t|~|\mathscr{F}_t^Y)
\xrightarrow{t\to\infty}0$ $\mathbf{P}^\nu$-a.s.
\end{lem}

\begin{proof}
Denote $|\varrho_\infty-\varrho_t|=u_t$ and 
$\mathbf{E}^\nu(|\varrho_\infty-\varrho_t|~|\mathscr{F}_t^Y)=v_t$,
and fix a constant $\varepsilon>0$.  Define the following stopping times:
$$
        \tau_1=\inf\{t>0:v_t>\varepsilon\},
        \qquad
        \sigma_1=\inf\{t>\tau_1:v_t<\varepsilon/2\},
$$
and for any $n\ge 2$
$$
        \tau_n=\inf\{t>\sigma_{n-1}:v_t>\varepsilon\},
        \qquad
        \sigma_n=\inf\{t>\tau_n:v_t<\varepsilon/2\}.
$$
By right-continuity of the sample paths, $v_{\tau_n}\ge\varepsilon$ on 
$\{\tau_n<\infty\}$.  But then
$$
        \varepsilon\,\mathbf{P}^\nu(\tau_n<\infty) = 
        \varepsilon\,\mathbf{P}^\nu(v_{\tau_n}
                \,I_{\tau_n<\infty}\ge\varepsilon)\le
        \mathbf{E}^\nu(v_{\tau_n}\,I_{\tau_n<\infty}),
$$
where we have used Chebyshev's inequality.  But as $v_t$ is the optional
projection of $u_t$ \cite[chapter VI, theorems 43 and 47]{DM}, we can write
$v_{\tau_n}\,I_{\tau_n<\infty}=\mathbf{E}^\nu(u_{\tau_n}
\,I_{\tau_n<\infty}|\mathscr{F}_{\tau_n}^Y)$ $\mathbf{P}^\nu$-a.s.  Hence,
in particular, $\varepsilon\,\mathbf{P}^\nu(\tau_n<\infty)
\le\mathbf{E}^\nu(u_{\tau_n}\,I_{\tau_n<\infty})$.

We now claim that $\tau_n\to\infty$ as $n\to\infty$ a.s.  To see this,
note that $\tau_n$ is nondecreasing, so it must converge either to
infinity or to a finite value.  But if it converges to a finite value,
then that sample path of $Y_t$ must have a discontinuity of the second
kind and hence cannot be c{\`a}dl{\`a}g.  Thus we can conclude that
$\tau_n\to\infty$ a.s., and hence $u_{\tau_n}\to 0$ a.s.\ by L{\'e}vy's
upward theorem (as $\varrho_t\to\varrho_\infty$ a.s.).  We would like to
show that $u_{\tau_n}\to 0$ in $L^1$, so that we can conclude that
$\mathbf{P}^\nu(\tau_n<\infty)\to 0$ as $n\to\infty$. To this end, note
that $u_{\tau_n}\to 0$ in $L^1$ is equivalent to
$\varrho_{\tau_n}\to\varrho_\infty$ in $L^1$.  But applying again
the optional projection property, we find that $\varrho_{\tau_n} =
\mathbf{E}^\nu(\varrho_\infty|\mathscr{F}^Y_{\tau_n})$
$\mathbf{P}^\nu$-a.s.  Hence the desired convergence follows from
L{\'e}vy's upward theorem.

We have established that $\mathbf{P}^\nu(\tau_n<\infty)\to 0$ as
$n\to\infty$. It follows directly that
$\mathbf{P}^\nu(\tau_n<\infty\mathrm{~for~all~}n)\le
\inf_n\mathbf{P}^\nu(\tau_n<\infty)=0$, so with unit probability either
$\limsup_{n\to\infty}v_t\le\varepsilon$, or
$\liminf_{n\to\infty}v_t\ge\varepsilon/2$.  But note that
$\|v_t\|_1\le\|u_t\|_1\to 0$ as $t\to\infty$, so $v_t\to 0$ in $L^1$.
Hence $\liminf_{n\to\infty}v_t\ge\varepsilon/2$ can only happen on a
null set, and we conclude that $\limsup_{n\to\infty}v_t\le\varepsilon$
a.s.  As this holds for any $\varepsilon>0$, the desired convergence is 
established.
\end{proof}

We are finally in a position to prove the main result.

\begin{thm}
\label{thm:stab}
Let $\mu\ll\nu$ and $f\in\mathcal{O}$.  Then
$$
        |\mathbf{E}^\mu(f(X_t)|\mathscr{F}_t^Y)-
        \mathbf{E}^\nu(f(X_t)|\mathscr{F}_t^Y)|\to 0
        \quad
        \mathbf{P}^\mu\mbox{-a.s.}
$$
\end{thm}

\begin{proof}
First, note that it suffices to prove the theorem for $f\in\mathcal{O}^0$.
After all, suppose we have established the result for $\mathcal{O}^0$.
By proposition \ref{prop:dense}, there is for $f\in\mathcal{O}$ a 
sequence $f_n\in\mathcal{O}^0$ such that $\|f-f_n\|\to 0$ as $n\to\infty$.  
Then
\begin{equation*}
\begin{split}
        &\limsup_{t\to\infty}
        |\mathbf{E}^\mu(f(X_t)|\mathscr{F}_t^Y)-
        \mathbf{E}^\nu(f(X_t)|\mathscr{F}_t^Y)| \\
        &\phantom{thi}
        \le \limsup_{t\to\infty}
        |\mathbf{E}^\mu(f(X_t)-f_n(X_t)|\mathscr{F}_t^Y)| 
        +\limsup_{t\to\infty}
        |\mathbf{E}^\nu(f_n(X_t)-f(X_t)|\mathscr{F}_t^Y)| \\
        &\phantom{the qu} +
        \limsup_{t\to\infty}
        |\mathbf{E}^\mu(f_n(X_t)|\mathscr{F}_t^Y)-
        \mathbf{E}^\nu(f_n(X_t)|\mathscr{F}_t^Y)| \\
        &\phantom{thi} \le
        2\,\|f-f_n\| + 
        \limsup_{t\to\infty}|\mathbf{E}^\mu(f_n(X_t)|\mathscr{F}_t^Y)-
        \mathbf{E}^\nu(f_n(X_t)|\mathscr{F}_t^Y)| \\
        &\phantom{thi}
        = 2\,\|f-f_n\|\quad\mathbf{P}^\mu\mbox{-a.s.}
\end{split}
\end{equation*}
But then the result follows for $f\in\mathcal{O}$ by letting $n\to\infty$.

We may thus assume that $f\in\mathcal{O}^0$, and by the linearity of the
conditional expectation we may assume without loss of
generality that $f$ is of the form
$$
        f(x) = \mathbf{E}_{(x,y)}(\xi),\qquad
        \xi = 
        f_1(Y_{t_1}-Y_0)f_2(Y_{t_2}-Y_0)\cdots f_n(Y_{t_n}-Y_0),
$$
for some $n<\infty$, $t_i\in D$ and bounded
continuous functions $f_i$.  By the Markov property, 
we find that $f(X_t) = \mathbf{E}^\nu(\xi\circ\theta_t|\mathscr{F}_t)$
$\mathbf{P}^\nu$-a.s.\ and that $f(X_t) = \mathbf{E}^\mu(
\xi\circ\theta_t|\mathscr{F}_t)$ $\mathbf{P}^\mu$-a.s., so we obtain
$$
        |\mathbf{E}^\mu(f(X_t)|\mathscr{F}_t^Y)-
        \mathbf{E}^\nu(f(X_t)|\mathscr{F}_t^Y)| =
        |\mathbf{E}^\mu(\xi\circ\theta_t|\mathscr{F}_t^Y)-
        \mathbf{E}^\nu(\xi\circ\theta_t|\mathscr{F}_t^Y)|.
$$
But as the $f_i$ are continuous, $\xi_t=\xi\circ\theta_t$ has 
c{\`a}dl{\`a}g sample paths, and clearly $\xi_t$ is
$\mathscr{F}^Y$-measurable for every $t$.  It remains to apply
lemma \ref{lem:cvgsimple}.
\end{proof}

An immediate consequence is that observability implies stability.

\begin{defn}
A filtering model is \textit{stable} if whenever $\mu\ll\nu$, 
$$
        |\mathbf{E}^\mu(f(X_t)|\mathscr{F}_t^Y)-
        \mathbf{E}^\nu(f(X_t)|\mathscr{F}_t^Y)|\to 0
        \quad
        \mathbf{P}^\mu\mbox{-a.s.}
        \quad
        \mbox{for all }
        f\in\mathcal{C}_b(\mathbb{S}).
$$
\end{defn}

\begin{cor}
\label{cor:stab}
If the filtering model is observable, then it is stable.
\end{cor}

\begin{proof}
This is immediate from the definition of observability.
\end{proof}

\begin{rem}
A word should be said at this point about the assumptions that the
signal process is a Markov process and that the observation process has
conditionally independent increments.  There is nothing essential in the
convergence proofs that depends on these properties, and indeed these
can safely be dropped (in fact, one may then choose the observation
state space $\mathbb{O}$ to be any locally compact Polish space).  In
this case, however, we could not guarantee that the space of observable
functions will contain only functions on $\mathbb{S}$; instead, we would
obtain $\mathcal{O}\subset\mathcal{C}_0(\mathbb{S}\times\mathbb{O})$ and
$\mathcal{N}\subset\mathcal{M}(\mathbb{S}\times\mathbb{O})$, and we
would have to consider convergence of conditional expectations of the
form $\mathbf{E}(f(X_t,Y_t)|\mathscr{F}_t^Y)$.  In other words, in this
case the initial measure on the observation process can play a nontrivial
role, which is not surprising.  The setting in which we have chosen to 
work---where the signal dynamics does not depend on the observations and
the observation noise is memoryless---is the natural setting where the 
initial measure on the observations decouples from the problem.
This allows us to concentrate on filtered estimates of the signal process,
which are the quantities which are of interest in the majority of applications.
\end{rem}

\begin{rem}
\label{rem:ultra}
Our notion of stability requires that $\mu\ll\nu$.  This is unavoidable
if we wish to define the filtered estimates as conditional expectations:
as $\mathbf{E}^\mu(f(X_t)|\mathscr{F}_t^Y)$ is only defined up to
$\mathbf{P}^\mu$-a.s.\ equivalence, the comparison of
$\mathbf{E}^\mu(f(X_t)|\mathscr{F}_t^Y)$ and
$\mathbf{E}^\nu(f(X_t)|\mathscr{F}_t^Y)$ for $\mu\not\ll\nu$ need not
make sense under any measure.  In many cases, however, there is a
natural version of the conditional expectations which may be defined
simultaneously with respect to all $\mathbf{P}^\nu$.  In this case, one
may ask whether the filter is \textit{strong stable}, i.e., whether
stability holds even for $\mu\not\ll\nu$.  This typically requires a
\textit{controllability} assumption in addition to observability (section
\ref{sec:ultrageneral}).  For the time being we are chiefly interested in 
observability, but we will return to the strong stability problem in section 
\ref{sec:ultra} in the setting of white noise type observations.
\end{rem}

\section{White noise type and counting observations}
\label{sec:obsmodels}

The purpose of this section is to investigate how two specific 
observation models that are extremely common in practice---white noise
and counting observations---fit into the general results developed 
in the previous subsections.

\subsection{White noise type observations}
\label{sec:whitenoise}

We consider the following setting:
$X_t$ is a Feller-Markov process, and $Y_t$ can be written in the form
$$
        Y_t = Y_0 + \int_0^t h(X_s)\,ds + K B_t,
$$
where $K$ is a non-random $p\times q$ matrix,
$h:\mathbb{S}\to\mathbb{O}$ is a continuous function, and $K B_t$ is a
$p$-dimensional Wiener process, with covariance matrix $KK^*$, which is
independent of $X_t$ and $Y_0$ (for any $\mathbf{P}_\nu$).  Note that
$K$ may be degenerate, in which case $v^*K B_t$ could be identically
zero for certain $v\in\mathbb{R}^p$.

\begin{lem}
The white noise type observation model satisfies the
conditionally independent increments property.
\end{lem}

\begin{proof}
Let $\xi$ be any bounded, $\mathscr{F}^X\vee\mathscr{G}^Y$-measurable
random variable.  Then $\xi$ is $\mathscr{F}^X\vee\sigma\{K
B_t:t>0\}$-measurable. We claim that for any $\mathscr{F}^X\vee\sigma\{K
B_t:t>0\}$-measurable random variable,
$\mathbf{E}_{(x,y)}(\xi)$ is independent of $y$.  To establish this, it
suffices to prove the claim for functions of the form $\xi_1\xi_2$ where
$\xi_1$ is $\mathscr{F}^X$-measurable and $\xi_2$ is $\sigma\{K
B_t:t>0\}$-measurable; the statement then follows by the monotone class
theorem.  But $\mathbf{E}_{(x,y)}(\xi_1\xi_2)=
\mathbf{E}_{(x,y)}(\xi_1)\mathbf{E}_{(x,y)}(\xi_2)$ by independence,
while $\mathbf{E}_{(x,y)}(\xi_1)$ only depends on $x$ (as $X_t$ is a
Markov process) and $\mathbf{E}_{(x,y)}(\xi_2)$ depends on neither $x$
or $y$ (as $B_t$ is a Wiener process for every $\mathbf{P}_{\nu}$).
\end{proof}

As the observations only depend on the signal through the observation
function $h$, a natural question is whether the observable and
nonobservable spaces depend on the noise covariance $KK^*$.  As one might
expect, this is not the case; for the purpose of observability, we may
simply take $K=0$. This is very convenient in computations, and shows
that observability is a structural property which does not depend on the
signal-to-noise ratio of the observations.

\begin{prop}
\label{prop:whiteredux}
For the white noise type observation model, 
$$
        \mathcal{N} = \{\alpha\mu_1-\alpha\mu_2\in\mathcal{M}(\mathbb{S}):
	\alpha\in\mathbb{R},~\mu_1,\mu_2\in\mathcal{P}(\mathbb{S}),
	~\mathbf{P}^{\mu_1}|_{\mathscr{G}^h}=
	\mathbf{P}^{\mu_2}|_{\mathscr{G}^h}\},
$$
where $\mathscr{G}^h=\sigma\{h(X_t):t\ge 0\}$.  
\end{prop}

To prove this statement, we will need the following simple lemma.

\begin{lem}
\label{lem:charf}
Let $(Z_1,\ldots,Z_n)$ and $(Z_1',\ldots,Z_n')$ be arbitrary random variables,
and let $(\xi_1,\ldots,\xi_n)$ be Gaussian random variables independent of all
$Z_i,Z_i'$.  Then 
$$
	(Z_1,\ldots,Z_n)\stackrel{\mathrm{law}}{=}
	(Z_1',\ldots,Z_n')
	\quad\mbox{iff}\quad
	(Z_1+\xi_1,\ldots,Z_n+\xi_n)
	\stackrel{\mathrm{law}}{=}(Z_1'+\xi_1,\ldots,Z_n'+\xi_n).
$$
\end{lem}

\begin{proof}
Recall that a probability measure on $\mathbb{R}^n$ is uniquely
determined by its characteristic function. 
Denote by $\chi_Z$, $\chi_{Z'}$, $\chi_{Z+\xi}$, $\chi_{Z'+\xi}$, and
$\chi_\xi$ the characteristic functions of $(Z_1,\ldots,Z_n)$,
$(Z_1',\ldots,Z_n')$, $(Z_1+\xi_1,\ldots,Z_n+\xi_n)$,
$(Z_1'+\xi_1,\ldots,Z_n'+\xi_n)$, and $(\xi_1,\ldots,\xi_n)$,
respectively.  Then, by independence, $\chi_{Z+\xi}= \chi_{Z}\chi_{\xi}$
and $\chi_{Z'+\xi}= \chi_{Z'}\chi_{\xi}$.  But as $\xi$ is a Gaussian
random vector, $\chi_{\xi}$ is invertible, so evidently
$\chi_{Z+\xi}=\chi_{Z'+\xi}$ iff $\chi_Z=\chi_{Z'}$.   This establishes
the claim.
\end{proof}

\begin{proof}[Proof of proposition \ref{prop:whiteredux}]
Recall that $\mu\in\mathcal{N}$ iff there exist probability measures
$\mu_1,\mu_2\in\mathcal{P}(\mathbb{S})$ and $\alpha>0$ such that
$\mu=\alpha\mu_1-\alpha\mu_2$ and $\mathbf{P}^{\mu_1}|_{\mathscr{F}^Y}=
\mathbf{P}^{\mu_2}|_{\mathscr{F}^Y}$. By \cite[theorem 16.6]{Bil} the
finite dimensional distributions form a separating class for probability
measures on $D([0,\infty[\mbox{};\mathbb{O})$, so $\mu\in\mathcal{N}$
iff $\mathbf{P}^{\mu_1}|_{\mathscr{G}}=
\mathbf{P}^{\mu_2}|_{\mathscr{G}}$ for all
$\mathscr{G}=\sigma\{Y_{t_1},\ldots,Y_{t_n}\}$ with $n<\infty$ and
$t_1,\ldots,t_n\in[0,\infty[\mbox{}$. But by lemma \ref{lem:charf}, this
is the case iff $\mathbf{P}^{\mu_1}|_{\mathscr{H}}=
\mathbf{P}^{\mu_2}|_{\mathscr{H}}$ for all
$\mathscr{H}=\sigma\{\int_0^{t_1}h(X_s)\,ds,\ldots,
\int_0^{t_n}h(X_s)\,ds\}$ with $n<\infty$ and
$t_1,\ldots,t_n\in[0,\infty[\mbox{}$.  As $\int_0^th(X_s)\,ds$ is
continuous (and in particular c{\`a}dl{\`a}g), so that the finite
dimensional distributions form a separating class also for this process,
and as $\sigma\{\int_0^th(X_s)\,ds:t>0\}=\mathscr{G}^h$, the result follows.  
\end{proof}

\subsection{Counting observations}

We now turn to the case of counting observations, for which almost
identical results hold.  In this setting, $X_t$ is again a Feller-Markov
process, and $Y_t=Y_0+N_t$ where $N_t$ is a Cox process
\cite[proposition 10.5]{Kal97} with intensity $\lambda_t=h(X_t)$,
conditionally independent of $Y_0$ given $\mathscr{F}^X$, for every
$\mathbf{P}_\mu$. By definition, this implies that for any $\mu$, under
a regular conditional probability
$\mathbf{P}_\mu(\,\cdot\,|\mathscr{F}^X)$ (which exists as our spaces
are Polish), $N_t^i$ are independent Poisson processes with intensities
$\lambda_t^i$ and the process $N_t$ is independent of $Y_0$.  Here the
observation function $h:\mathbb{S}\to\mathbb{O}$ is a continuous
nonnegative function.

\begin{lem}
The counting observation model satisfies the
conditionally independent increments property.
\end{lem}

\begin{proof}
Let $\xi$ be any bounded, $\mathscr{F}^X\vee\mathscr{G}^Y$-measurable
random variable. Then by our assumptions, under a regular conditional
probability $\mathbf{P}_{(x,y)}(\,\cdot\,|\mathscr{F}^X)$, the law of
$\xi$ only depends on the sample paths of $X_t$ and is thus independent
of $y$. In particular, this means that
$\mathbf{E}_{(x,y)}(\xi|\mathscr{F}^X)$ does not depend on $y$.  But
then $\mathbf{E}_{(x,y)}(\xi)=
\mathbf{E}_{(x,y)}(\mathbf{E}_{(x,y)}(\xi|\mathscr{F}^X))$ can not
depend on $y$, as $X_t$ is a Markov process in its own right and as
$\mathbf{E}_{(x,y)}(\xi|\mathscr{F}^X)$ is an $\mathscr{F}^X$-measurable
random variable. 
\end{proof}

An analog of proposition \ref{prop:whiteredux} also holds.

\begin{prop}
\label{prop:poissredux}
The conclusion of proposition \ref{prop:whiteredux} holds identically for
the counting observation model.
\end{prop}

\begin{proof}
This follows directly from \cite[lemma 10.8]{Kal97}.
\end{proof}

\subsection{A simple sufficient condition}

Let us mention a useful consequence of these results, which leads
to a particularly simple sufficient condition for observability.

\begin{lem}
\label{lem:onetoone}
For the white noise type and counting observations models, it is always
the case that $f\circ h\in\mathcal{O}$ for any measurable function 
$f:\mathbb{O}\to\mathbb{R}$ such that $f\circ h\in\mathcal{C}_b(\mathbb{S})$.
In particular, if $h$ is one-to-one then we may conclude that
$\mathcal{O}=\mathcal{C}_b(\mathbb{S})$ (i.e., the signal-observation model 
is observable).
\end{lem}

\begin{proof}
Let $f\circ h\in\mathcal{C}_b(\mathbb{S})$, and choose any
$\mu\in\mathcal{N}$. Then $\mu=\alpha\mu_1-\alpha\mu_2$, where
$\mathbf{P}^{\mu_1}|_{\mathscr{G}^h}=
\mathbf{P}^{\mu_2}|_{\mathscr{G}^h}$.  Thus in particular
$\mathbf{P}^{\mu_1}|_{\sigma\{h(X_0)\}}=
\mathbf{P}^{\mu_2}|_{\sigma\{h(X_0)\}}$. But $f(h(X_0))$ is
$\sigma\{h(X_0)\}$-measurable, so that evidently
$$
        \int f(h(x))\,\mu_1(dx) = \int f(h(x))\,\mu_2(dx).
$$
As this holds for any $\mu\in\mathcal{N}$, we find that
$f\circ h\in\mathcal{O}$.
\end{proof}

In other words, ``nice'' functions of the observation function are
always observable, regardless of any further properties of the model. 

\begin{rem}
For the special case where $f$ is chosen to be the identity, the
stability of the observation function (in a slightly different sense)
was found in \cite[theorem 3.1]{COC99} under much weaker conditions. 
However, the latter result cannot be used to conclude the stability of
the filter, even in the case when $h$ is one-to-one. 
\end{rem}

\section{Finite state signals}
\label{sec:fss}

The simplest nonlinear filtering model is one where the signal state
space consists of a finite number of points.  Such models are
of particular theoretical and practical interest as the filtered
estimates can be finite dimensionally computed.  On the other hand, this
model shares many of the features of more general models and thus serves
as a convenient prototype.  In this section, we will use the results of
the previous section to obtain an essentially complete characterization
of the stability of such filters in the case of nondegenerate white
noise type observations.

Throughout this section, $X_t$ is a Markov process on the finite state
space $\mathbb{S}=\{a_1,\ldots,a_d\}$ with transition intensities
matrix $\Lambda=(\lambda_{ij})_{1\le i,j\le d}$.  The observations
process $Y_t$ is taken to be one-dimensional and of the form 
$$
        Y_t = Y_0 + \int_0^t h(X_s)\,ds+\kappa\,W_t,\qquad
        Y_0 = 0,
$$
where $h:\mathbb{S}\to\mathbb{R}$.  The restriction to one dimension is
for notational convenience only; all the results extend directly to
observations in $\mathbb{R}^p$.

\begin{rem}
$Y_0$ is included as a reminder that this is a Markov
observation model. As usual, we will enforce $Y_0=0$ by working with the
measures $\mathbf{P}^\mu$.
\end{rem}

In subsection \ref{sec:finiteobs}, we elaborate on the structure of the
spaces $\mathcal{N}$ and $\mathcal{O}$ in this setting (the results of
this section hold identically for the case of counting observations).
Subsection \ref{sec:finitedet} (see also section \ref{sec:finiteult}) is
devoted to the complete characterization of the stability of the filter. 
Here we make essential use of the white noise type observations, and it
is moreover crucial that the observations are assumed to be
\textit{nondegenerate} $\kappa>0$.  The reason for this is that we will
invoke results that hold only in this setting; see remark
\ref{rem:degfs} below for further details.

Finally, a word on notation.  We use the following notation for the
filter: $\pi_t^\mu(f) = \mathbf{E}^{\mu}(f(X_t)|
\mathscr{F}_t^Y)$.  When $\pi_t^\mu$ is used as a vector, this is
implied in the sense that $(\pi_t^\mu)^i = \pi_t^\mu(I_{\{a_i\}})$.  We
will interchangeably treat functions on $\mathbb{S}$ as vectors in
$\mathbb{R}^d$ in the obvious way ($v^i=v(a_i)$), whenever this is
convenient.  The transpose of a vector or matrix is denoted as $v^*$ or
$M^*$.

\subsection{Observability}
\label{sec:finiteobs}

For the particular case of a finite state signal and $\kappa=0$, the
notion of observability has been investigated in the context of
identifiability and lumpability of hidden Markov models
\cite{obs1,obs2,obs3}, though chiefly in discrete time. In this
subsection, we briefly develop the necessary results in our setting.

Let $\mu,\nu\in\mathcal{P}(\mathbb{S})$.  To determine the nonobservable
space $\mathcal{N}$, by proposition \ref{prop:whiteredux}, we need to
find all $\mu,\nu$ with $\mathbf{P}^{\mu}|_{\mathscr{G}^h}=
\mathbf{P}^{\nu}|_{\mathscr{G}^h}$. But as $h(X_t)$ is a c{\`a}dl{\`a}g
process, it suffices to verify that the finite dimensional distributions
of $h(X_t)$ are the same under $\mathbf{P}^{\mu}$ and $\mathbf{P}^{\nu}$
\cite[theorem 16.6]{Bil}. We thus begin by computing these
distributions.

\begin{lem}
\label{lem:backcond}
Let $\mathbb{H}=h(\mathbb{S})=\{b_1,\ldots,b_r\}$, $r\le d$ be the set
of possible observation values.  Define the $d\times d$ projection
matrices $H_{b_k}$ such that $(H_{b_k})_{i,j}=1$ whenever $i=j$ and
$h(a_i)=b_k$, and zero otherwise.  Then under
$\mathbf{P}^{\mu}$, the finite dimensional
distributions of $h(X_t)$ have the form
\begin{multline*}
        \mathbf{P}^{\mu}(h(X_0)=n_0,h(X_{t_1})=n_1,\ldots,
                h(X_{t_k})=n_k) = \\
        \mu^*H_{n_0}e^{\Lambda t_1}H_{n_1}e^{\Lambda(t_2-t_1)}H_{n_2}
        \cdots e^{\Lambda(t_k-t_{k-1})}H_{n_k}\mathbf{1},
\end{multline*}
where $n_i\in\mathbb{H}$ and $\mathbf{1}\in\mathbb{R}^d$ is the vector of 
ones.
\end{lem}

\begin{proof}
The result follows from
\begin{multline*}
        \mathbf{P}^{\mu}(X_0=m_0,X_{t_1}=m_1,\ldots,
                X_{t_k}=m_k) = \\
        \mu_{m_0}[e^{\Lambda t_1}]_{m_0m_1}
        [e^{\Lambda(t_2-t_1)}]_{m_1m_2}\cdots
        [e^{\Lambda(t_k-t_{k-1})}]_{m_{k-1}m_k}.
\end{multline*}
by summing $m_j$ over the set $\{a_i\in\mathbb{S}:h(a_i)=n_j\}$.
\end{proof}
We immediately conclude the following.
\begin{cor}
The observable and nonobservable spaces satisfy
\begin{equation*}
\begin{split}
        \mathcal{O} &= \mathrm{span}\left\{
        H_{n_0}e^{\Lambda\delta_1}H_{n_1}e^{\Lambda\delta_2}H_{n_2}
        \cdots e^{\Lambda\delta_k}H_{n_k}\mathbf{1}:
        k\ge 0,~\delta_i>0,~n_i\in\mathbb{H}\right\}, \\
        \mathcal{N} &= \mathcal{O}^\perp =
        \{v\in\mathbb{R}^d:
        v^*x=0\mbox{ for all }x\in\mathcal{O}\}.
\end{split}
\end{equation*}
The model is observable if and only if $\dim\mathcal{O}=d$.
\end{cor}

The following simplification is useful in computations.

\begin{lem}
\label{lem:obsshort}
The observable space can be characterized as follows:
$$
        \mathcal{O}=\mathrm{span}\left\{
        H_{n_0}\Lambda H_{n_1}\Lambda
        \cdots \Lambda H_{n_k}\mathbf{1}:
        k\ge 0,~n_i\in\mathbb{H}\right\}.
$$
\end{lem}

\begin{proof}
Note that any vector of the form ($p_i\le d-1$)
$$
        H_{n_0}\Lambda^{p_1}H_{n_1}\Lambda^{p_2}H_{n_2}
        \cdots \Lambda^{p_k}H_{n_k}\mathbf{1}
$$
can be obtained from a vector of the form
$$
        H_{n_0}e^{\Lambda\delta_1}H_{n_1}e^{\Lambda\delta_2}H_{n_2}
        \cdots e^{\Lambda\delta_k}H_{n_k}\mathbf{1}
$$
by taking derivatives with respect to $\delta_i$, and in particular the
former is the limit of elements of $\mathcal{O}$.  But $\mathcal{O}$ is
closed as it is a finite dimensional linear space, so the span of the 
former is contained in $\mathcal{O}$.  To prove the converse 
inclusion, it suffices to expand the matrix exponential in a power series 
and apply the Cayley-Hamilton theorem.  Finally, note that $H_{b_i}$ sum
to the identity matrix, so we can reduce to the case where $p_i=1$ for all $i$.
\end{proof}

We will need, in particular, the following important consequence.

\begin{cor}
$\mathcal{O}$ is invariant under $\Lambda$ and $H_{b_i}$:
$\Lambda\mathcal{O}\subset\mathcal{O}$,
$H_{b_i}\mathcal{O}\subset\mathcal{O}$.
Similarly $\mathcal{N}$ is invariant under $\Lambda^*$ and 
$H_{b_i}$: $\Lambda^*\mathcal{N}\subset\mathcal{N}$,
$H_{b_i}\mathcal{N}\subset\mathcal{N}$.
\end{cor}

\begin{proof}
Immediate from the previous lemma.
\end{proof}

\begin{rem}
Given the previous corollary, it is not surprising that $\mathcal{O}$
can in fact be characterized by its invariance property.  To this end,
denote by 
$$
	\mathcal{O}_h=\{f\circ
	h:\forall\,f:\mathbb{H}\to\mathbb{R}\}=
	\mathrm{span}\{H_{b_{i}}\mathbf{1}:b_i\in\mathbb{H}\}.
$$ 
Then $\mathcal{O}$ is the smallest subspace of $\mathbb{R}^d$ that
contains $\mathcal{O}_h$ and is invariant under $\Lambda$ and all
$H_{b_i}$.  Indeed, let us call this smallest subspace $\mathcal{O}'$.
Clearly $\mathcal{O}'\subset\mathcal{O}$, as $\mathcal{O}$ contains
$\mathcal{O}_h$ and is invariant under $\Lambda$ and $H_{b_i}$. On the
other hand, every element of $\mathcal{O}$ can be generated from
elements in $\mathcal{O}_h$ by a finite number of multiplications by
$\Lambda$ and $H_{b_i}$ and linear combinations.  Hence
$\mathcal{O}\subset\mathcal{O}'$.

To verify observability, we could proceed as follows.  Denote
$$
        \mathcal{Z}_1 = \mathcal{O}_h,\qquad
        \mathcal{Z}_n = \mathcal{Z}_{n-1}+
                \Lambda\mathcal{Z}_{n-1}+
                H_{b_1}\mathcal{Z}_{n-1}+\cdots+
                H_{b_r}\mathcal{Z}_{n-1},\quad
        n>1,
$$
where the sum of two linear spaces denotes their linear span.  It is
evident that every element of $\mathcal{O}$ will be in $\mathcal{Z}_n$
for some $n$.  Moreover, if $\mathcal{Z}_n=\mathcal{Z}_{n+1}$ for some
$n=m$, then it is true for all $n>m$, and in particular $\mathcal{Z}_m=
\mathcal{O}$.  Finally, we claim that this will always be the case
for some $m<d$.  Indeed, the dimension of $\mathcal{Z}_n$ can not shrink
with increasing $n$, but it can not grow larger than $d$ as we are working in
$\mathbb{R}^d$.  As $\mathcal{O}_h$ contains at least the constants,
the procedure must complete in at most $d-1$ steps.  This idea is classical, 
see, e.g., \cite[section 3.2.2]{geom}, and could be implemented, e.g., 
by starting with the natural basis $\{H_{b_i}\mathbf{1}\}$ of $\mathcal{O}_h$ 
and applying the Gram-Schmidt procedure at every iteration $n$ to obtain a 
basis for $\mathcal{Z}_n$.  
\end{rem}

\begin{rem}
In an early paper on filter stability, Delyon and Zeitouni \cite{delyon}
impose a condition (A2) which, by the previous remark, is seen to be
sufficient (but not necessary) for observability.  In addition, they
assume ergodicity of the signal process.  Though their condition (A2)
was later shown to be superfluous \cite[theorem 4.1]{BCL} in the
nondegenerate case $\kappa>0$, Delyon and Zeitouni show through a
counterexample that when their condition (A2) is not satisfied, the
filter may lose its stability as $\kappa\to 0$.  That this can not
happen when condition (A2) is satisfied is to be expected as, by corollary
\ref{cor:stab} above, observability implies filter stability without any
nondegeneracy or ergodicity assumptions. It does not appear, however, that
our results can be related to the methods used in \cite{delyon}, nor do our 
results give any information on the rate of convergence (exponential 
convergence is proved in \cite{delyon}).
\end{rem}

\begin{rem}
\label{rem:wononedim}
Denote by $C=[H_{b_1}\mathbf{1},\ldots,H_{b_r}\mathbf{1}]$ the $d\times
r$ matrix whose columns are indicator functions on level sets of $h$.
A sufficient (but not necessary) condition for observability is that
$\mathrm{rank}([C~\Lambda C~\Lambda^2C~\cdots~\Lambda^{d-1}C])=d$,
which is the classical observability test for linear systems.
This corresponds to considering only the one dimensional distributions
of $h(X_t)$, rather than all finite dimensional distributions.
\end{rem}

\subsection{A complete characterization of filter stability when
$\kappa>0$}
\label{sec:finitedet}

Corollary \ref{cor:stab} and the results of the previous section
show that
\begin{cor}
If $\dim\mathcal{O}=d$, then the filter is stable.
\end{cor}
The converse, however, is not true.  In the nondegenerate case
$\kappa>0$, it was shown by Baxendale, Chigansky and Liptser \cite{BCL}
that ergodicity of the signal is a sufficient condition for stability of
the filter, regardless of the observation structure.  It is not
difficult to find an example of a filtering model that is not
observable, but has an ergodic signal (e.g., choose any ergodic signal
and set $h=0$; another example is the one in \cite{delyon}).

The goal of this section is to find a necessary and sufficient condition
for filter stability \textit{in the nondegenerate case $\kappa>0$},
which we will assume throughout.  This is done by combining our results
above with the results from \cite{BCL}.  To gain some intuition, recall
that $|\pi_t^\mu(f)-\pi_t^\mu(f)|\to 0$, or, equivalently,
$|(\pi_t^\mu)^*f-(\pi_t^\nu)^*f|\to 0$, for any $f\in\mathcal{O}$. This
implies that as $t\to\infty$, the signed measure $\pi_t^\mu-\pi_t^\nu$
converges to the nonobservable space $\mathcal{N}$. To ensure stability,
we would like to find a condition under which the space $\mathcal{N}$
converges to zero under the dynamics of the filter.

One plausible condition is to require that the signal itself ``forgets''
perturbations in $\mathcal{N}$, i.e., that $|\mathbf{P}^{\mu}(X_t=a_i)-
\mathbf{P}^{\nu}(X_t=a_i)|\to 0$ as $t\to\infty$ for all
$a_i\in\mathbb{S}$ whenever $\mu-\nu\in\mathcal{N}$. In this way, we
obtain the natural counterpart of the notion of \textit{detectability}
in linear systems theory.  We will show that this condition is indeed
necessary and sufficient for stability of the filter, provided that $\kappa>0$.

Before turning to the proof, let us make precise what we are going to show.
Recall that $\Lambda^*\mathcal{N}\subset\mathcal{N}$; hence it makes sense
to speak of the restriction $\Lambda^*|_{\mathcal{N}}$.

\begin{lem}
\label{lem:hurwitz} 
Denote $(p_t^\mu)^i=\mathbf{P}^{\mu}(X_t=a_i)$,
and suppose that we have $\dim\mathcal{N}>0$.
Then the following are equivalent statements.
\begin{enumerate}
\item $|p_t^\mu-p_t^\nu|\to 0$ as $t\to\infty$ whenever
$\mu-\nu\in\mathcal{N}$.
\item $\Lambda^*|_{\mathcal{N}}$ is Hurwitz (its eigenvalues
have strictly negative real parts).
\item $\Lambda^*|_{\mathcal{N}}$ has full rank.
\end{enumerate}
Here $|v|$ denotes the $\ell_1$-norm of the vector $v$.
\end{lem}

\begin{proof}
The Kolmogorov forward equation states that
$$
        \frac{d}{dt}\,p_t^\mu = \Lambda^*p_t^\mu,\qquad
        p_0^\mu=\mu.
$$
Hence, in particular,
$$
        \frac{d}{dt}(p_t^\mu-p_t^\nu) = 
        \Lambda^*(p_t^\mu-p_t^\nu) =
        \Lambda^*|_{\mathcal{N}}(p_t^\mu-p_t^\nu),\qquad
        p_0^\mu-p_0^\nu=\mu-\nu\in\mathcal{N},
$$
where the second equality follows immediately from the fact that this
equation leaves $\mathcal{N}$ invariant.  It is well known from linear
systems theory that the solution of this equation decays to zero as
$t\to\infty$ for every initial condition $\mu-\nu\in\mathcal{N}$ if and
only if $\Lambda^*|_{\mathcal{N}}$ is Hurwitz.  The fact that
$\Lambda^*|_\mathcal{N}$ is Hurwitz if and only if it is of full rank
follows from the fact that any nonzero eigenvalue of $\Lambda^*$ has
strictly negative real part \cite[pages 52--53]{Bar55}.
\end{proof}

Our previous discussion now motivates the following definition.

\begin{defn}
The signal-observation model is called \textit{detectable} if it is
either observable or any of the equivalent conditions of lemma
\ref{lem:hurwitz} hold.
\end{defn}

The goal of this section is to prove the following theorem.

\begin{thm}
\label{thm:detfin}
Suppose that $\kappa>0$.  Then
$|\pi_t^\mu-\pi_t^\nu|\xrightarrow{t\to\infty}0$ $\mathbf{P}^{\mu}$-a.s.\ 
whenever $\mu\ll\nu$ if and only if the signal-observation model is detectable.
\end{thm}

\begin{rem}
\label{rem:degfs}
The situation for $\kappa=0$ appears to be more complicated, and the
theorem does not hold in this case.  A counterexample can be found in
\cite[section 3]{BCL} (see also \cite{delyon}), which discusses a model
that is certainly detectable, but the filter is not stable when
$\kappa=0$ due to a sort of ``geometric obstruction''.  Problems of this
sort, in somewhat different setting, date back to the work of Kaijser 
\cite{kaijser}, and some recent progress on that problem can be found in
\cite{reeds}.  A complete understanding of this case is still lacking, however.
\end{rem}

\subsubsection{Necessity}

To prove theorem \ref{thm:detfin}, we begin by showing that
detectability is a necessary condition for the stability of the filter.

\begin{lem}
\label{lem:necessary}
Suppose that $|\pi_t^\mu-\pi_t^\nu|\xrightarrow{t\to\infty}0$
$\mathbf{P}^{\mu}$-a.s.\ for any $\mu\ll\nu$.  Then
the signal-observation model is detectable.
\end{lem}

\begin{proof}
Let $\mu-\nu\in\mathcal{N}$; then $\mathbf{P}^{\mu}$ and
$\mathbf{P}^{\nu}$ are identical on $\mathscr{F}^Y$.  As $\pi_t^\nu(f)$
is $\mathscr{F}^Y$-measurable, $\mathbf{E}^{\mu}(\pi_t^\nu(f))=
\mathbf{E}^{\nu}(\pi_t^\nu(f))$. Thus
$$
        \mathbf{E}^{\mu}|\pi_t^\mu(f)-\pi_t^\nu(f)|\ge
        |\mathbf{E}^{\mu}(\pi_t^\mu(f)-\pi_t^\nu(f))|
        =
        |\mathbf{E}^{\mu}(f(X_t))-
                \mathbf{E}^{\nu}(f(X_t))|\quad
        \mbox{for any }f.
$$
Now suppose the model is not detectable, i.e., there exists 
$v\in\mathcal{N}$ so that $|p_t^\mu-p_t^\nu|\not\to 0$ when $\mu-\nu\propto v$.
Choose $\nu=(v^+ + v^-)/2v^+(\mathbb{S})$ and $\mu=v^+/v^+(\mathbb{S})$;
then $\mu,\nu\in\mathcal{P}(\mathbb{S})$, $\mu\ll\nu$ and $\mu-\nu\propto v$.
As $p_t^\mu-p_t^\nu$ does not converge to zero, there must
exist a function $f\in\mathcal{C}_b(\mathbb{S})$ and a sequence of times 
$t_n\nearrow\infty$ such that
$$
        |\mathbf{E}^{\mu}(f(X_{t_n}))
        -\mathbf{E}^{\nu}(f(X_{t_n}))|
        \xrightarrow{n\to\infty} \alpha>0.
$$
But if $|\pi_t^\mu-\pi_t^\nu|\xrightarrow{t\to\infty}0$
$\mathbf{P}^{\mu}$-a.s., then using dominated convergence
$$
        |\mathbf{E}^{\mu}(f(X_{t_n}))
        -\mathbf{E}^{\nu}(f(X_{t_n}))|
        \le \mathbf{E}^{\mu}
                |\pi_{t_n}^\mu(f)-\pi_{t_n}^\nu(f)|
        \xrightarrow{n\to\infty}0.
$$
Hence we have a contradiction, and the proof is complete.
\end{proof}

\begin{rem}
The previous proof does not use at all the fact that $\mathbb{S}$ is a finite
set or that the observations are of the white noise type.  Indeed, let us call
a general model \textit{detectable} if $\mu-\nu\in\mathcal{N}$ implies that 
$|\mathbf{E}^\mu(f(X_t))-\mathbf{E}^\nu(f(X_t))|\to 0$ as
$t\to\infty$ for any $f\in\mathcal{C}_b(\mathbb{S})$.  Then precisely the
same proof shows that detectability is a \textit{necessary} condition for 
the stability of the filter (it is not even necessary to assume that 
$\mathbb{S}$ is compact).  The difficult part is to establish that 
detectability is a \textit{sufficient} condition for stability of the filter, 
and this is what we will do below for finite state signals with nondegenerate
white noise type observations.
\end{rem}

\subsubsection{Sufficiency: no transient states}

We now proceed to prove that detectability is also a sufficient
condition for stability when $\kappa>0$.  \textit{Throughout this and the
following subsection we always assume that the signal-observation model is
detectable}.

In the proofs, we make use of the partition of the state space
$\mathbb{S}$ into $M<\infty$ ergodic classes $\mathbb{U}_i$,
$i=1,\ldots,M$ and a transient class $\mathbb{T}$, so that $\mathbb{S}$
is the disjoint union of these sets.  Any Markov chain can be uniquely
decomposed in this way. 

\begin{lem}
If $\Lambda F=0$ for a function $F$, then $F\in\mathcal{O}$.
\end{lem}

\begin{proof}
$\Lambda F=0$ is equivalent to $e^{\Lambda t}F=F$ for all $t\ge 0$.  
In particular, this implies that $(p_t^\mu)^*F=\mu^*F$ for all $t\ge 0$.  
Now suppose that $F\not\in\mathcal{O}$; then there exist $\mu,\nu$ such that
$\mu-\nu\in\mathcal{N}$ and $|\mu^*F-\nu^*F|=\alpha>0$.  In particular,
we find that $|(p_t^\mu)^*F-(p_t^\nu)^*F|=\alpha$ for all $t\ge 0$.
But detectability implies that $|p_t^\mu-p_t^\nu|\to 0$ as $t\to\infty$
for $\mu-\nu\in\mathcal{N}$, so we have a contradiction.
\end{proof}

\begin{cor}
Suppose $\mathbb{T}=\varnothing$.  Then 
$I_{\mathbb{U}_i}\in\mathcal{O}$
for any $i=1,\ldots,M$.
\end{cor}

\begin{proof}
It is easily seen that $\Lambda I_{\mathbb{U}_i}=0$ when there are no
transient states.  Hence the statement follows from the previous lemma.
\end{proof}

Suppose that $\mathbb{T}=\varnothing$. The essential consequence of
detectablity is that as $t\to\infty$, we will be able to determine
precisely in which of the ergodic classes
$\mathbb{U}_1,\ldots,\mathbb{U}_M$ the signal started at $t=0$. 
Following the logic of \cite{BCL}, this will cause the filter to be
stable when combined with the fact that the filter is stable for ergodic
signals.  We will deal with the transient states separately in the next
subsection, and assume for now that there are no such states (or,
equivalently, that we work with initial densities that are supported
on the ergodic classes only).

\begin{lem}
\label{lem:proofofsmallmeas}
$\mathbf{E}^{\nu}(I_{\mathbb{U}_i}(X_0)|\mathscr{F}_t^Y)
\xrightarrow{t\to\infty}I_{\mathbb{U}_i}(X_0)$
$\mathbf{P}^{\nu}$-a.s.,
provided that there are no transient states
$\mathbb{T}=\varnothing$.
\end{lem}

\begin{proof}
For any $j$ such that $\nu(\mathbb{U}_j)>0$, denote by 
$\nu_j=I_{\mathbb{U}_j}\nu/\nu(\mathbb{U}_j)$.
Then $\nu_j\ll\nu$ and $\pi_t^{\nu_j}(I_{\mathbb{U}_i})=\delta_{ij}$.
But by theorem \ref{thm:stab}
$|\pi_t^{\nu_j}(I_{\mathbb{U}_i})-\pi_t^\nu(I_{\mathbb{U}_i})|\to 0$
as $t\to\infty$ $\mathbf{P}^{\nu_j}$-a.s.,
as $I_{\mathbb{U}_i}\in\mathcal{O}$.  In other words,
$\mathbf{P}^{\nu}(\pi_t^\nu(I_{\mathbb{U}_i})\not\to
\delta_{ij}\mbox{ and }X_0\in\mathbb{U}_j)=0$, 
so $\pi_t^\nu(I_{\mathbb{U}_i})\to\delta_{ij}$ on
$\{\omega:X_0\in\mathbb{U}_j\}$, modulo a
$\mathbf{P}^{\nu}$-null set.  Finally, note that
$I_{\mathbb{U}_i}(X_t)=I_{\mathbb{U}_i}(X_0)$ $\mathbf{P}^\nu$-a.s., as the 
ergodic classes do not communicate.
\end{proof}

We can now prove sufficiency for the special case $\mathbb{T}=\varnothing$.

\begin{lem}
\label{lem:detectnotrans}
Suppose $\mathbb{T}=\varnothing$ and $\mu\ll\nu$.  Then
$|\pi_t^\mu(f)-\pi_t^\nu(f)|\to 0$ $\mathbf{P}^{\mu}$-a.s.
for all $f\in\mathcal{C}_b(\mathbb{S})$.
\end{lem}

\begin{proof}
By the Bayes formula, we find that $\mathbf{P}^{\mu}$-a.s.
$$
        \mathbf{E}^{\mu}(f(X_t)|\mathscr{F}_t^Y)=
        \sum_{j=1}^M
        \mathbf{E}^{\mu_j}(f(X_t)|\mathscr{F}_t^Y)\,
        \mathbf{P}^{\mu}(X_0\in\mathbb{U}_j|
                \mathscr{F}_t^Y).
$$
The same equation holds with $\mu,\mu_j$ replaced by $\nu,\nu_j$.
The result now follows easily from the previous lemma and the fact that
$|\mathbf{E}^{\mu_j}(f(X_t)|\mathscr{F}_t^Y)-
\mathbf{E}^{\nu_j}(f(X_t)|\mathscr{F}_t^Y)|\to 0$ by
\cite[theorem 4.1]{BCL} (as $X_t$ is supported entirely in the ergodic
class $\mathbb{U}_j$ under the initial measures $\mu_j,\nu_j$).
\end{proof}

\subsubsection{Sufficiency: general case}

We now consider the general case with $\mathbb{T}\ne\varnothing$.  Let
us begin by showing that the transient states themselves decay as
$t\to\infty$.

\begin{lem}
\label{lem:notransientsplease}
$\pi_t^\nu(I_\mathbb{T})\to 0$ as $t\to\infty$
$\mathbf{P}^{\nu}$-a.s.
\end{lem}

\begin{proof}
Note that $I_{\mathbb{T}}(X_t)\to 0$ 
$\mathbf{P}^{\nu}$-a.s., as the transient states must
decay eventually into one of the ergodic classes.  Now write
$$
        \pi_t^\nu(I_{\mathbb{T}}) 
        = \mathbf{E}^{\nu}\left(\left.
                I_{\mathbb{T}}(X_t)\right|\mathscr{F}_t^Y\right) \le
        \mathbf{E}^{\nu}\left(\left.
                \sup_{u\ge n}I_{\mathbb{T}}(X_u)\right|\mathscr{F}_t^Y\right)   
$$
for all $t\ge n$ $\mathbf{P}^{\nu}$-a.s.\ (using the 
c{\`a}dl{\`a}g paths to eliminate the time dependence of the null set).
Hence
$$
        \limsup_{t\to\infty}\mathbf{E}^{\nu}\left(\left.
                I_{\mathbb{T}}(X_t)\right|\mathscr{F}_t^Y\right) \le
        \mathbf{E}^{\nu}\left(\left.
                \sup_{u\ge n}I_{\mathbb{T}}(X_u)\right|\mathscr{F}^Y
                \right).
$$
The claim follows by letting $n\to\infty$ using dominated convergence.
\end{proof}

Evidently, as $t\to\infty$, the conditional measures $\pi_t^\nu$ and
$\pi_t^\mu$ converge to measures that are supported on the ergodic
classes $\mathbb{U}=\mathbb{S}\backslash\mathbb{T}=
\bigcup_{i=1}^M\mathbb{U}_i$.  On the other hand, if we start with
$\mu\ll\nu$ which are already supported on $\mathbb{U}$, then
$|\pi_t^\mu-\pi_t^\nu|\to 0$ by lemma \ref{lem:detectnotrans}.  This
strongly suggests that we should have $|\pi_t^\mu-\pi_t^\nu|\to 0$ for
any $\mu\ll\nu$.  Our goal is to prove this assertion.

\begin{lem}
\label{lem:bigtransest}
Suppose that $\mu\ll\nu$ and that $\mu$ is supported on $\mathbb{U}$.
Then
$$
        \mathbf{E}^{\mu}\left(
                \limsup_{t\to\infty}|\pi_t^{\mu}(f)-\pi_t^{\nu}(f)|
        \right)\le \mathrm{osc}(f)~\|d\mu/d\nu\|_\infty~\nu(\mathbb{T}),
$$
where $\mathrm{osc}(f) = \max(f)-\min(f)$.
\end{lem}

\begin{proof}
Let us write $\nu_{\mathbb{U}}=I_{\mathbb{U}}\nu/\nu(\mathbb{U})$ and
$\nu_{\mathbb{T}}=I_{\mathbb{T}}\nu/\nu(\mathbb{T})$.
By the Bayes formula, we find that $\mathbf{P}^{\nu}$-a.s.
$$
        \pi_t^\nu(f) =
        \pi_t^{\nu_{\mathbb{U}}}(f)\,
        \mathbf{P}^{\nu}(X_0\in\mathbb{U}|\mathscr{F}_t^Y)
        +
        \pi_t^{\nu_{\mathbb{T}}}(f)\,
        \mathbf{P}^{\nu}(X_0\in\mathbb{T}|\mathscr{F}_t^Y).
$$
It follows directly that
$\mathbf{P}^{\nu_{\mathbb{U}}}$-a.s.
$$
        |\pi_t^{\nu_{\mathbb{U}}}(f)-\pi_t^\nu(f)|
        = 
        |\pi_t^{\nu_{\mathbb{U}}}(f)-\pi_t^{\nu_{\mathbb{T}}}(f)|\,
        \mathbf{P}^{\nu}(X_0\in\mathbb{T}|\mathscr{F}_t^Y),
$$
so that in particular
$$
        \limsup_{t\to\infty}|\pi_t^{\nu_{\mathbb{U}}}(f)-\pi_t^\nu(f)|
        \le
        \mathrm{osc}(f)\,
        \mathbf{P}^{\nu}(X_0\in\mathbb{T}
                |\mathscr{F}^Y).
$$
We thus compute
$$
        \mathbf{E}^{\mu}(
        \mathbf{P}^{\nu}(X_0\in\mathbb{T}
                |\mathscr{F}^Y)) 
         \le
        \|d\mu/d\nu\|_\infty~
        \mathbf{E}^{\nu}(
        \mathbf{P}^{\nu}(X_0\in\mathbb{T}
                |\mathscr{F}^Y)) 
         = \|d\mu/d\nu\|_\infty~\nu(\mathbb{T}).
$$
The claim now follows from lemma \ref{lem:detectnotrans} and
$$
        \limsup_{t\to\infty}|\pi_t^\mu(f)-\pi_t^\nu(f)|
        \le
        \limsup_{t\to\infty}|\pi_t^\mu(f)-\pi_t^{\nu_{\mathbb{U}}}(f)| +
        \limsup_{t\to\infty}|\pi_t^{\nu_{\mathbb{U}}}(f)-\pi_t^\nu(f)|,
$$
using the fact that $\mu$ and $\nu_{\mathbb{U}}$ are both supported on 
$\mathbb{U}$.
\end{proof}

To establish that the right-hand side in the expression in this lemma
can be chosen to be zero, we will use the Markov property
of the filter. 

\begin{lem}
\label{lem:wonhamfeller}
For $\mu\ll\nu$, the pair $(\pi_t^\mu,\pi_t^\nu)$ is a
Feller-Markov process under $\mathbf{P}^{\mu}$.
\end{lem}

\begin{proof}
Recall that $dB_t^\mu = \kappa^{-1}(dY_t - \pi_t^\mu(h)\,dt)$, the
innovations process, is an $\mathscr{F}_t^Y$-Wiener process under
$\mathbf{P}^{\mu}$, and that we thus have
\begin{equation*}
\begin{split}
        & d\pi_t^\mu = \Lambda^*\pi_t^\mu\,dt +
                \kappa^{-1}\,(H-h^*\pi_t^\mu)\,\pi_t^\mu\,dB_t^\mu,\\
        & d\pi_t^\nu = \Lambda^*\pi_t^\nu\,dt +
                \kappa^{-1}\,(H-h^*\pi_t^\nu)\,\pi_t^\nu\,
        (dB_t^\mu + \kappa^{-1}\,h^*\pi_t^\mu\,dt -
                \kappa^{-1}\,h^*\pi_t^\nu\,dt),
\end{split}
\end{equation*}
where $H=\mathrm{diag}(h)$ and $(\pi_0^\mu,\pi_0^\nu)=(\mu,\nu)$.
For these facts, see, e.g., \cite{LS01a}.  Being the solution of a stochastic
differential equation with Lipschitz coefficients (the coefficients are 
bounded in the double simplex $\Delta^d\times\Delta^d$, and the first
exit time from the simplex is infinite), it is well known that there
is a unique strong solution which satisfies the Markov and Feller
properties.
\end{proof}

A particular consequence of this lemma is the following.  Consider
the pair of $\Delta^d$-valued stochastic differential equations
\begin{equation*}
\begin{split}
        & d\pi_t = \Lambda^*\pi_t\,dt +
                \kappa^{-1}\,(H-h^*\pi_t)\,\pi_t\,dW_t,\\
        & d\bar\pi_t = \Lambda^*\bar\pi_t\,dt +
                \kappa^{-1}\,(H-h^*\bar\pi_t)\,\bar\pi_t\,
        (dW_t + \kappa^{-1}\,h^*\pi_t\,dt -
                \kappa^{-1}\,h^*\bar\pi_t\,dt),
\end{split}
\end{equation*}
where $W_t$ is a standard Wiener process.  The solutions of this
stochastic differential equation can be realized on the canonical path
space $\tilde\Omega = D([0,\infty[\mbox{};\Delta^d)\times
D([0,\infty[\mbox{};\Delta^d)$ such that $\pi_t(u,v)=u(t)$ and
$\bar\pi_t(u,v)=v(t)$ are the canonical processes, and with a family of
measures $\mathbb{P}_{(\mu,\nu)}$ under which $(\pi_t,\bar\pi_t)$ solve
the stochastic differential equation above for the initial condition
$(\pi_0,\bar\pi_0)=(\mu,\nu)$.  We can subsequently introduce the
natural filtration $\mathscr{E}_t=\sigma\{(\pi_s,\bar\pi_s): s\le t\}$,
augmented as usual with respect to the family $\mathbb{P}_{(\mu,\nu)}$,
and the canonical shift $\tilde\theta_t(u,v)(s)= (u(s+t),v(s+t))$, such
that the process $(\pi_t,\bar\pi_t)$ satisfies the usual Markov property
with respect to the filtration $\mathscr{E}_t$ and the family
$\mathbb{P}_{(\mu,\nu)}$.  From the proof of the previous lemma, it
follows that for any $\mu\ll\nu$, the law of the process
$(\pi_t,\bar\pi_t)$ under $\mathbb{P}_{(\mu,\nu)}$ coincides with the
law of the process $(\pi_t^\mu,\pi_t^\nu)$ under $\mathbf{P}^\mu$. In
particular, our previous results can be applied to the process
$(\pi_t,\bar\pi_t)$, and to establish stability it suffices to
demonstrate the corresponding property for the latter. 

\begin{rem} 
The construction of $(\pi_t,\bar\pi_t)$ on its own path space is
certainly not necessary, but helps alleviate some notational confusion. 
In particular, we will be using the Markov property of the filter,
whereas our previous notation is geared at the Markov property of the
signal-observation pair. 
\end{rem}

Combined with lemma \ref{lem:bigtransest}, we can now establish the following.

\begin{lem}
\label{lem:prefinalwon}
Suppose that $\mu\ll\nu$ and that $\mu$ is supported on $\mathbb{U}$.
Then it follows that $|\pi_t^\mu(f)-\pi_t^{\nu}(f)|\to 0$
$\mathbf{P}^{\mu}$-a.s.\ for any $f\in\mathcal{C}_b(\mathbb{S})$.
\end{lem}

\begin{proof}
Using the Markov property, we can write
$$
	\mathbb{E}_{(\mu,\nu)}\left(\left.
		\limsup_{t\to\infty}|\pi_t(f)-\bar\pi_t(f)|
	~\right|\mathscr{E}_s\right) =
	\mathbb{E}_{(\pi_s,\bar\pi_s)}\left(
		\limsup_{t\to\infty}|\pi_t(f)-\bar\pi_t(f)|
	\right),
$$
where we have used the fact that the random variable
$\limsup_{t\to\infty}|\pi_t(f)-\bar\pi_t(f)|$ is invariant
under the shift $\tilde\theta_s$.  By \cite[lemma 2.1]{COC99}
we find that $\pi_s\ll\bar\pi_s$ $\mathbb{P}_{(\mu,\nu)}$-a.s.\ whenever
$\mu\ll\nu$, whereas clearly $\pi_s$ is $\mathbb{P}_{(\mu,\nu)}$-a.s.\
supported on $\mathbb{U}$ whenever $\mu$ is supported on $\mathbb{U}$.
Hence we can invoke lemma \ref{lem:bigtransest}, and we find that
$$
	\mathbb{E}_{(\mu,\nu)}\left(\left.
		\limsup_{t\to\infty}|\pi_t(f)-\bar\pi_t(f)|
	~\right|\mathscr{E}_s\right) \le
        \mathrm{osc}(f)~
        \|d\pi_s/d\bar\pi_s\|_\infty~\bar\pi_s(I_\mathbb{T}),
$$
where $\|d\pi_s/d\bar\pi_s\|_\infty = \max_i(d\pi_s/d\bar\pi_s)^i$.
In particular, this implies that
$$
        \mathbf{E}^{\mu}\left(\left.
                \limsup_{t\to\infty}|\pi_t^\mu(f)-\pi_t^{\nu}(f)|
        ~\right|\mathscr{F}_s^Y
        \right)\le \mathrm{osc}(f)~
        \|d\pi_s^\mu/d\pi_s^\nu\|_\infty~\pi_s^\nu(I_\mathbb{T}).
$$
Now note that (see, e.g., \cite[lemma 2.1]{COC99})
$$
        \left\|\frac{d\pi_s^\mu}{d\pi_s^\nu}\right\|_\infty =
        \max_i\,
        \frac{\mathbf{E}^{\nu}(
                \frac{d\mu}{d\nu}(X_0)|\mathscr{F}_s^Y,X_s=a_i)}{
        \mathbf{E}^{\nu}(
                \frac{d\mu}{d\nu}(X_0)|\mathscr{F}_s^Y)}
        \le \frac{1}{\varrho_s}
        \left\|\frac{d\mu}{d\nu}\right\|_\infty.
$$
Hence, letting $s\to\infty$ and using
lemma \ref{lem:notransientsplease}, we find that
$$
        \limsup_{t\to\infty}|\pi_t^\mu(f)-\pi_t^{\nu}(f)|=0
        \quad
        \mbox{on}\quad\{\omega:\varrho_\infty(\omega)>0\}\backslash N,
$$
where $\mathbf{P}^{\mu}(N)=0$.  But
$\mathbf{P}^{\mu}(\varrho_\infty>0)=1$, so the result follows.
\end{proof}

We can now finally complete the proof.  It is important to remember that
we have assumed detectability throughout this subsection.

\begin{prop}
\label{prop:finallywon}
Suppose the signal-observation model is detectable.
If $\mu\ll\nu$, then $|\pi_t^\mu(f)-\pi_t^{\nu}(f)|\to 0$
$\mathbf{P}^{\mu}$-a.s.\ for any $f\in\mathcal{C}_b(\mathbb{S})$.
\end{prop}

\begin{proof}
By the previous lemma, we find that
$|\pi_t^{\mu_{\mathbb{U}}}(f)-\pi_t^\mu(f)|\to 0$
and $|\pi_t^{\nu_{\mathbb{U}}}(f)-\pi_t^\nu(f)|\to 0$ 
$\mathbf{P}^{\mu_{\mathbb{U}}}$-a.s.  Hence, using
the triangle inequality, $|\pi_t^\mu(f)-\pi_t^\nu(f)|\to 0$
$\mathbf{P}^{\mu_{\mathbb{U}}}$-a.s.  But this implies,
as in the proof of lemma \ref{lem:proofofsmallmeas}, that
$|\pi_t^\mu(f)-\pi_t^\nu(f)|\to 0$ on $\{\omega:X_0\in\mathbb{U}\}$,
modulo a $\mathbf{P}^{\mu}$-null set.  In particular, we can then estimate
$$
        \mathbf{E}^{\mu}\left(
        \limsup_{t\to\infty}|\pi_t^\mu(f)-\pi_t^\nu(f)|\right)
        \le \mathrm{osc}(f)~\mu(\mathbb{T}).
$$
To proceed, we apply the Markov property as in the previous proof.
This yields
$$
        \mathbf{E}^{\mu}\left(\left.
        \limsup_{t\to\infty}|\pi_t^\mu(f)-\pi_t^\nu(f)|
        ~\right|\mathscr{F}_s^Y\right)
        \le \mathrm{osc}(f)~\pi_s^\mu(I_\mathbb{T}).
$$
The result follows by letting $s\to\infty$ and using
lemma \ref{lem:notransientsplease}.
\end{proof}

\section{Strong stability for nondegenerate white noise type observations}
\label{sec:ultra}

\subsection{Strong stability}

Up to this point, we have always assumed that the initial measures of
interest are absolutely continuous $\mu\ll\nu$.  In this section we
consider the case when $\mu\not\ll\nu$. As explained in remark
\ref{rem:ultra}, the filter stability problem is in general not even
well defined for such initial measures, and the characterization of
strong stability (the stability of the filter for arbitrary initial
conditions) requires choosing a particular version of the conditional
expectations.  In the case of nondegenerate white noise type
observations, however, there is a natural choice of version, viz.\ the
one provided by the Kallianpur-Striebel formula, whose construction we
now briefly recall (see, e.g., \cite[section 7.9]{LS01a}).

We consider the generic white noise type observation model
$$
        Y_t = Y_0 + \int_0^t h(X_s)\,ds + K B_t
$$
of section \ref{sec:whitenoise}.  Beside the assumptions of section
\ref{sec:whitenoise}, however, we additionally assume
\textit{nondegeneracy} of the observation process, i.e., we assume the
$K$ is an invertible matrix.  The importance of this requirement stems
from the fact that it allows us to use Girsanov's theorem to remove the
dependence of the observations on the signal process, a fact which will
be exploited shortly.

Denote by $\mathbf{Q}^\mu$ the measure on the space of signal sample
paths $\Omega^X$ such that the canonical process $X_t(x)=x(t)$ has the
law of the signal process with initial distribution $\mu$, i.e.,
$\mathbf{Q}^\mu$ is the marginal of $\mathbf{P}^\mu$ on $\Omega^X$. 
Moreover, we denote by $\mathbf{W}$ the Wiener measure on $\Omega^Y$
with covariance $KK^*$ (i.e., $Y_t(y)=y(t)$ is a Wiener process with
covariance $KK^*$ under $\mathbf{W}$).  Then, by Girsanov's theorem,
$$
	Z_t =
	\frac{d\,\mathbf{P}^\mu|_{\mathscr{F}_t}
		}{d\,(\mathbf{Q}^\mu\times\mathbf{W})|_{\mathscr{F}_t}} =
	\exp\left(
		\int_0^t (KK^*)^{-1}h(X_s)\cdot dY_s 
		- \frac{1}{2}\int_0^t \|K^{-1}h(X_s)\|^2\,ds
	\right).
$$
Thus, using the Bayes formula, we obtain the following characterization of
the filter:
$$
	\pi_t^\mu(f)(x,y) = \mathbf{E}^\mu(f(X_t)|\mathscr{F}_t^Y)(x,y) =
	\frac{
		\int_{\Omega^X} Z_t(\tilde x,y)\,f(\tilde x(t))\,
			\mathbf{Q}^\mu(d\tilde x)
	}{
		\int_{\Omega^X} Z_t(\tilde x,y)\,\mathbf{Q}^\mu(d\tilde x)
	}\quad
	\mathbf{P}^\mu\mbox{-a.s.}
$$
This is the Kallianpur-Striebel formula.  Note, however, that the
expression on the right hand side depends only on the observation sample
paths $\Omega^Y$ in the time interval $[0,t]$, and is well defined not
only $\mathbf{P}^\mu$-a.s.\ but in fact for $\mathbf{W}$-a.e.\
$y\in\Omega^Y$.  Moreover, it is easily seen from Girsanov's theorem
that the observation marginals satisfy $\mathbf{P}^\mu|_{\mathscr{F}_T^Y}\sim
\mathbf{W}|_{\mathscr{F}_T^Y}$ for any $T<\infty$, regardless of the initial 
measure $\mu$.  Hence the Kallianpur-Striebel formula defines a version of 
$\pi_t^\mu$ which is $\mathbf{P}^\nu$-a.s.\ uniquely defined for any $\nu$ 
(even when $\mu$ and $\nu$ are not absolutely continuous).  \textit{In the 
remainder of this section, $\pi_t^\mu$ will always imply this particular 
version.}

Having now chosen a version of the filter that is well defined under any
measure $\mathbf{P}^\mu$, strong stability can be meaningfully defined.

\begin{defn}
The filtering model is \textit{strong stable} if
for any $\mu,\nu,\gamma\in\mathcal{P}(\mathbb{S})$,
$$
	|\pi_t^\mu(f)-\pi_t^\nu(f)|\xrightarrow{t\to\infty}0
	\quad	\mathbf{P}^\gamma\mbox{-a.s.}\quad
	\mbox{for all }f\in\mathcal{C}_b(\mathbb{S}).
$$
\end{defn}

\subsection{A complete characterization in the finite state case}
\label{sec:finiteult}

In the finite state setting, the condition $\mu\ll\nu$ is not really
restrictive in practice.  Indeed, if we wish to ensure that the filter
is asymptotically insensitive to its initial condition, it suffices to
initialize the filter with a (possibly incorrect) initial distribution
$\nu$ which charges every point in the state space, e.g., the uniform
distribution on $\mathbb{S}$.  As any measure on $\mathbb{S}$ is
absolutely continuous with respect to such a measure, the convergence of
the thus initialized filter to the optimal one is ensured, regardless of
the true initial distribution $\mu$, provided that the model is detectable and
$\kappa>0$.

Nonetheless the strong stability property is of interest, and can be
characterized completely as we did for the stability property. Somewhat
surprisingly, the observation structure no longer plays a role in this
setting.

\begin{thm}
\label{thm:ultdetfin}
Suppose that $\kappa>0$.  Then $|\pi_t^\mu-\pi_t^\nu|\xrightarrow{t\to\infty}0$
$\mathbf{P}^\gamma$-a.s.\ for any $\mu,\nu,\gamma$ if and only if the signal 
process has only one ergodic class.
\end{thm}

Before proceding with the proof of the theorem, let us make the following
important remark.

\begin{rem}
The stochastic differential equations used in the proof of lemma 
\ref{lem:wonhamfeller} can be obtained directly from the Kallianpur-Striebel
formula, and hence define the version of $(\pi_t^\mu,\pi_t^\nu)$ which we
use in this section.  In particular, this implies that the arguments based
on the Markov property of $(\pi_t^\mu,\pi_t^\nu)$ continue to hold in the
current setting and the condition $\mu\ll\nu$ of lemma \ref{lem:wonhamfeller}
is no longer required.
\end{rem}

Let us first prove the necessity part of the theorem.  We assume
throughout that $\kappa>0$ and that $\pi_t^\mu,\pi_t^\nu$ are chosen to
be the Kallianpur-Striebel versions.

\begin{lem}
Suppose that $|\pi_t^\mu-\pi_t^\nu|\xrightarrow{t\to\infty}0$
$\mathbf{P}^\gamma$-a.s.\ for any $\mu,\nu,\gamma$.  Then the signal
process has only one ergodic class.
\end{lem}

\begin{proof}
Suppose that the signal process has two ergodic classes $\mathbb{U}_1$
and $\mathbb{U}_2$.  Choose $\mu$ to be any distribution that is
supported on $\mathbb{U}_1$, and $\nu$ to be any distribution that is
supported on $\mathbb{U}_2$.  Then it is easily verified that
$\pi_t^\mu(I_{\mathbb{U}_1})=1$ $\mathbf{W}$-a.s.\ for all times $t$,
while $\pi_t^\nu(I_{\mathbb{U}_1})=0$ $\mathbf{W}$-a.s.\ for all times $t$.
Hence we have a contradiction.
\end{proof}

We now proceed to prove sufficiency.  First, note that it suffices to prove
that $|\pi_t^\mu-\pi_t^\nu|\to 0$ $\mathbf{P}^\mu$-a.s.\ for any $\mu,\nu$.
Indeed, it then follows that
$$
	|\pi_t^\mu-\pi_t^\nu| \le
	|\pi_t^\mu-\pi_t^\gamma| + |\pi_t^\nu-\pi_t^\gamma|
	\xrightarrow{t\to\infty}0
	\quad\mathbf{P}^\gamma\mbox{-a.s.}
$$
for any $\mu,\nu,\gamma$ by the triangle inequality.  As before, it is easier
to first consider the case with no transient states $\mathbb{T}=\varnothing$.

\begin{lem}
Suppose the signal process is ergodic (in particular $\mathbb{T}=\varnothing$).
Then $|\pi_t^\mu-\pi_t^\nu|\to 0$ $\mathbf{P}^\mu$-a.s.\ for any $\mu,\nu$.
\end{lem}

\begin{proof}
This is precisely \cite[theorem 4.1]{BCL}.
\end{proof}

Moreover, we need the following lemma.

\begin{lem}
Suppose the signal process has only one ergodic class $\mathbb{U}$.  Then
$\pi_t^\mu(a_i)>0$ a.s.\ for all $a_i\in\mathbb{U}$ and $t>0$, regardless
of $\mu$.
\end{lem}

\begin{proof}
By \cite[eq.\ (7.205)]{LS01a}, we have $\pi_t^\mu(a_i)>0$ a.s.\ if and
only if $\mathbf{P}^\mu(X_t=a_i)>0$.  But in the absence of multiple ergodic
classes, it must be the case that $\mathbf{P}^\mu(X_t=a_i)>0$ for any 
$a_i\in\mathbb{U}$ as soon as $t>0$, regardless of $\mu$.
\end{proof}

The transient states can now be eliminated precisely as in proposition
\ref{prop:finallywon}, completing the proof.

\begin{lem}
Suppose the signal process has only one ergodic class. Then we have
$|\pi_t^\mu-\pi_t^\nu|\to 0$ $\mathbf{P}^\mu$-a.s.\ for any $\mu,\nu$.
\end{lem}

\begin{proof}
First, note that as in the proof of lemma \ref{lem:prefinalwon},
we can write
$$
	\mathbb{E}_{(\mu,\nu)}\left(\left.
		\limsup_{t\to\infty}|\pi_t(f)-\bar\pi_t(f)|
	~\right|\mathscr{E}_s\right) =
	\mathbb{E}_{(\pi_s,\bar\pi_s)}\left(
		\limsup_{t\to\infty}|\pi_t(f)-\bar\pi_t(f)|
	\right),
$$
where we have used the Markov property and the fact that the random
variable $\limsup_{t\to\infty}|\pi_t(f)-\bar\pi_t(f)|$ is invariant
under the shift $\tilde\theta_s$.  But by the previous lemma, we find
that $(\pi_t)_{\mathbb{U}}\sim(\bar\pi_t)_{\mathbb{U}}$ for any $t>0$. 
Thus we may assume, without loss of generality, that
$\mu_{\mathbb{U}}\sim\nu_{\mathbb{U}}$ in the following.

By lemma \ref{lem:prefinalwon}, we find that
$|\pi_t^{\mu_{\mathbb{U}}}(f)-\pi_t^\mu(f)|\to 0$
$\mathbf{P}^{\mu_{\mathbb{U}}}$-a.s., and also that
$|\pi_t^{\nu_{\mathbb{U}}}(f)-\pi_t^\nu(f)|\to 0$
$\mathbf{P}^{\nu_{\mathbb{U}}}$-a.s.  But from the assumption
$\mu_{\mathbb{U}}\sim\nu_{\mathbb{U}}$, it follows that
$\mathbf{P}^{\mu_{\mathbb{U}}}\sim\mathbf{P}^{\nu_{\mathbb{U}}}$. Hence,
using the triangle inequality, $|\pi_t^\mu(f)-\pi_t^\nu(f)|\to 0$
$\mathbf{P}^{\mu_{\mathbb{U}}}$-a.s.  But this implies, as in the proof
of lemma \ref{lem:proofofsmallmeas}, that
$|\pi_t^\mu(f)-\pi_t^\nu(f)|\to 0$ on $\{\omega:X_0\in\mathbb{U}\}$,
modulo a $\mathbf{P}^{\mu}$-null set.  In particular, we can then
estimate
$$
        \mathbf{E}^{\mu}\left(
        \limsup_{t\to\infty}|\pi_t^\mu(f)-\pi_t^\nu(f)|\right)
        \le \mathrm{osc}(f)~\mu(\mathbb{T}).
$$
To proceed, we apply the Markov property.  This yields
$$
        \mathbf{E}^{\mu}\left(\left.
        \limsup_{t\to\infty}|\pi_t^\mu(f)-\pi_t^\nu(f)|
        ~\right|\mathscr{F}_s^Y\right)
        \le \mathrm{osc}(f)~\pi_s^\mu(I_\mathbb{T}).
$$
The result follows by letting $s\to\infty$ and using
lemma \ref{lem:notransientsplease}.
\end{proof}

\subsection{A general criterion}
\label{sec:ultrageneral}

In this section, we will give a sufficient condition for strong stability
for a general signal process (with nondegenerate white noise
observations as above). We will need the following definitions.

\begin{defn}
\label{def:regular}
Let $R_t$ be a Markov semigroup on a Polish space with associated
transition probabilities $p_t(x,A)$.  Then $R_t$ is called 
\begin{itemize}
\item \textit{strong Feller} if $R_tf$ is continuous for every bounded 
measurable $f$;
\item \textit{irreducible} if for any nonempty open set $A$, we have 
$p_t(x,A)>0$ for all $x$;
\item \textit{regular} if $p_t(x,\cdot)\sim p_t(y,\cdot)$ for every
$x,y$ and $t>0$.
\end{itemize}
\end{defn}

A well known result of Has'minski{\u\i} \cite{Hasm60} states the following.

\begin{lem}
Any irreducible strong Feller semigroup is regular.
\end{lem}

The following theorem and its immediate corollary are the main results
of this section.  Recall that, by assumption, the signal is a Markov
process in its own right.

\begin{thm}
\label{thm:genultra}
If the signal process is regular, then stability implies strong stability.
\end{thm}

\begin{cor}
\label{cor:genultra}
Regularity of the signal and observability imply strong stability of the
filter.  In particular, the filter is strong stable if the signal-observation
model is observable and the signal process is irreducible and strong Feller.
\end{cor}

For the proof of the theorem, we need the following counterpart of
lemma \ref{lem:wonhamfeller}.

\begin{lem}
The pair $(\pi_t^\mu,\pi_t^\nu)$ is a
Feller-Markov process under $\mathbf{P}^{\mu}$.
\end{lem}

\begin{proof}
This follows as in the proof of \cite[theorem 2.3]{Kun71}.
The details are omitted.
\end{proof}

This implies that as in the finite state case, we can construct the
filter on its canonical path space $\tilde\Omega =
D([0,\infty[\mbox{};\mathcal{P}(\mathbb{S}))\times
D([0,\infty[\mbox{};\mathcal{P}(\mathbb{S}))$ (here
$\mathcal{P}(\mathbb{S})$ is endowed with the topology of weak
convergence, which turns it into a compact Polish space).  To be
precise, denote by $\mathbb{P}_{(\mu,\nu)}$ the probability measure on
$\tilde\Omega$ under which the canonical processes $\pi_t(u,v)=u(t)$ and
$\bar\pi_t(u,v)=v(t)$ have the same law as do $\pi_t^\mu$ and
$\pi_t^\nu$ under $\mathbf{P}^\mu$. As before we introduce the natural
filtration $\mathscr{E}_t=\sigma\{(\pi_s,\bar\pi_s): s\le t\}$,
augmented with respect to the family $\mathbb{P}_{(\mu,\nu)}$, and the
canonical shift $\tilde\theta_t(u,v)(s)= (u(s+t),v(s+t))$. It then
follows that the process $(\pi_t,\bar\pi_t)$ satisfies the usual Markov
property with respect to the filtration $\mathscr{E}_t$ and the family
$\mathbb{P}_{(\mu,\nu)}$.

The strategy for proving theorem \ref{thm:genultra} is now
straightforward. What we will show is that if the signal process is
regular, then $\pi_t^\mu\sim\pi_t^\nu$ a.s.\ for any $t>0$, regardless of
$\mu$ and $\nu$.  Using the Markov property of the filter, the strong stability
problem then reduces to the ordinary stability problem.

\begin{proof}[Proof of theorem \ref{thm:genultra}]
Denote by $\mathbf{Q}_t^\mu$ the law of $X_t$ under $\mathbf{P}^\mu$.
From regularity, it follows that $\mathbf{Q}_t^\mu\sim\mathbf{Q}_t^\nu$
for any $\mu,\nu$ and $t>0$.  But from the Kallianpur-Striebel formula,
it follows directly that $\pi_t^\mu\sim\mathbf{Q}_t^\mu$ a.s.\ with
$$
	\frac{d\pi_t^\mu}{d\mathbf{Q}_t^\mu}(z) =
	\frac{
		\int_{\Omega^X} Z_t(\tilde x,\cdot)\,
			\mathbf{Q}^\mu(d\tilde x|\tilde x(t)=z)
	}{
		\int_{\Omega^X} Z_t(\tilde x,\cdot)\,\mathbf{Q}^\mu(d\tilde x)
	}.
$$
Hence evidently $\pi_t^\mu\sim\pi_t^\nu$ a.s.\ for any $\mu,\nu$ and $t>0$.
Now note that for $f\in\mathcal{C}_b(\mathbb{S})$, 
$$
	\mathbb{E}_{(\mu,\nu)}\left(\left.
		\limsup_{t\to\infty}|\pi_t(f)-\bar\pi_t(f)|
	~\right|\mathscr{E}_s\right) =
	\mathbb{E}_{(\pi_s,\bar\pi_s)}\left(
		\limsup_{t\to\infty}|\pi_t(f)-\bar\pi_t(f)|
	\right),
$$
where we have used the Markov property and the fact that the random
variable $\limsup_{t\to\infty}|\pi_t(f)-\bar\pi_t(f)|$ is invariant
under the shift $\tilde\theta_s$.  But we have just established that
$\pi_t\sim\bar\pi_t$ a.s., and thus the right hand side of this expression
vanishes a.s.\ due to the fact that the filter is already assumed to be
stable.  Thus $\limsup_{t\to\infty}|\pi_t(f)-\bar\pi_t(f)|=0$ a.s., and 
the claim is established.
\end{proof}

The regularity of the signal process is closely related to the classical
notion of controllability.  Suppose that $\mathbb{S}$ is a compact
connected $C^\infty$-manifold, and that the signal process $X_t$ is the
solution of the Stratonovich stochastic differential equation
$$
	dX_t = F(X_t)\,dt + G(X_t)\circ dW_t,\qquad X_0\in\mathbb{S},
$$
where $F$ and $G$ are $C^\infty$-vector fields on $\mathbb{S}$.  Then
$X_t$ is a Markov process as usual with transition probabilities 
$p_t(x,A)$.  Consider also the associated control system
$$
	\frac{d}{dt}\,\Xi_t = F(\Xi_t) + G(\Xi_t)\,u(t),\qquad
		\Xi_0\in\mathbb{S},
$$
where $u(t)$ is the control input.  We denote by
$A_t(x)\subset\mathbb{S}$ the set of points $\Xi_t$ which are reachable
from $\Xi_0=x$ by the application of a piecewise smooth control signal
$u$, and call the signal \textit{controllable} if $A_t(x)=\mathbb{S}$
for every $x\in\mathbb{S}$ and $t>0$.  It follows from the
Stroock-Varadhan support theorem that controllability is a sufficient
condition for irreducibility of the signal \cite{Kun78}. Moreover, the
controllability assumption additionally implies hypoellipticity of the
diffusion $X_t$ (see the remark on \cite[page 175]{Kun78}), which gives
rise to the strong Feller property. 

Thus evidently a sufficient condition for strong stability of the filter,
for a diffusion signal on a compact manifold with white noise type
observations, is that the signal is controllable and the filtering model
is observable.  This mirrors precisely the well known
controllability-observability criterion for the stability of the Kalman
filter \cite{BJ68,OP96}.  Indeed, it is not difficult to verify that the
linear filtering model is observable in the sense of this paper
precisely when the well known observability rank condition is satisfied,
while the linear signal is controllable precisely when the
controllability rank condition is satisfied (though, unfortunately, the
Kalman filter does not fit into the current setting as its state space
is not compact).

\begin{rem}
Regularity of the signal process is certainly not a necessary condition
for strong stability.  In the finite state setting, for example,
regularity occurs only when there is a single ergodic class and there
are no transient states. We have seen, however, that strong stability
still holds true in the presence of transient states.  The latter
situation is analogous to the stabilizability criterion for the
stability of the Kalman filter \cite{OP96}.  One might hope that also
stabilizability and detectability have natural counterparts in the
general setting, but we will not pursue this here.

On the other hand, we remark that stabilizability and detectability are
generally considered together in the stability theory for the Kalman
filter, while the results of this paper indicate that these conditions
play rather separate roles. In particular, it is to be expected that
detectability is a sufficient condition for the stability of the Kalman
filter even in the absence of stabilizability, provided that the initial
distributions are absolutely continuous $\mu\ll\nu$.  That this is
indeed the case (under slightly stronger conditions) is shown in the
appendix.\footnote{
	However, the method used in the appendix to prove stability
	of the Kalman filter is unrelated to the techniques developed
	in the body of this paper.}
\end{rem}

Both the Kalman filter and the finite state case give rise to conditions
for observability and controllability which are easily computed
explictly in terms of matrices.  For general diffusions, the matter
appears to be much more complicated.  To establish controllability one
may employ certain Lie-algebraic computations, as detailed in
\cite{Kun78}.  The question of observability for signals on a non-finite
state space does not appear to have been studied at all in the
literature.

A slightly stronger condition than observability, however, is closely
related to the classical observability problem for (deterministic)
infinite-dimensional linear systems. Suppose that we have white noise
type or counting observations, so that observability is determined by
$\mathscr{G}^h=\sigma\{h(X_t):t\ge 0\}$. Rather than require every
$\mu\ne\nu$ to merely give rise to different
$\mathbf{P}^\mu|_{\mathscr{G}^h}\ne\mathbf{P}^\nu|_{\mathscr{G}^h}$, we
could ask whether $\mu\ne\nu$ implies that
$\mathbf{P}^\mu|_{\sigma\{h(X_t)\}}\ne\mathbf{P}^\nu|_{\sigma\{h(X_t)\}}$
for some $t\ge 0$.  The latter is clearly a sufficient condition for
observability, where only the one-dimensional distributions of the
process $h(X_t)$ are taken into account (compare with remark 
\ref{rem:wononedim}).

Now denote by $R_t$ the Markov semigroup of the signal process.  Then
there exists a dual semigroup $R_t^*$, which acts on the space of
measures $\mathcal{M}(\mathbb{S})$, such that $R_t^*\mu$ is the law of
$X_t$ under $\mathbf{P}^\mu$.  Moreover, let us define the projection
map $\Pi:\mathcal{M}(\mathbb{S})\to\mathcal{M}(h(\mathbb{S}))$ such that
$\Pi:\mu\mapsto \mu\circ h^{-1}$.  Then we can consider $\mathsf{X}_t =
R_t^*\mu$ as defining the dynamics of an infinite-dimensional linear
system with $\mathsf{X}_0=\mu$ and with infinite-dimensional linear
observations $\mathsf{Y}_t=\Pi\mathsf{X}_t$.  The classical
observability problem associated with this infinite-dimensional linear
model characterizes precisely when $\mu\ne\nu$ implies that
$\mathbf{P}^\mu|_{\sigma\{h(X_t)\}}\ne\mathbf{P}^\nu|_{\sigma\{h(X_t)\}}$
for some $t\ge 0$.  A detailed treatment of observability problems of
this type can be found in \cite{Trig76}.

\section{The non-compact case}
\label{sec:non-compact}

In the preceding sections we have considered exclusively the case where
the signal state space $\mathbb{S}$ is compact, so that
$\mathcal{C}(\mathbb{S})=\mathcal{C}_b(\mathbb{S})=\mathcal{C}_0(\mathbb{S})$.
When $\mathbb{S}$ is not compact, as in the common setting where
$\mathbb{S}=\mathbb{R}^n$, for example, one could try to extend the
proofs to show the stability of functions in
$\mathcal{C}_0(\mathbb{S})$. The latter space of functions is the
obvious choice from the point of view of our techniques, as
$\mathcal{C}_0(\mathbb{S})^*=\mathcal{M}(\mathbb{S})$ even when
$\mathbb{S}$ is only locally compact.  However, from a practical point
of view the stability of functions in $\mathcal{C}_0(\mathbb{S})$ is too
restrictive; indeed, if the signal is transient then the filtered
estimate for any such function will be stable, but this fact is of
little interest (as the filtered estimates of functions that vanish at
infinity yield no information on a transient signal as $t\to\infty$). 
Instead, one should consider the larger class of continuous bounded
functions $\mathcal{C}_b(\mathbb{S})$ or of all continuous
functions $\mathcal{C}(\mathbb{S})$.

Unfortunately, the techniques which we have developed in the previous
sections do not extend directly to this setting. The problem is, of
course, that the dual of $\mathcal{C}_b(\mathbb{S})$ (with respect to
the uniform topology) is no longer $\mathcal{M}(\mathbb{S})$ when
$\mathbb{S}$ is not compact; rather, $\mathcal{C}_b(\mathbb{S})^*$ can
be characterized as $\mathcal{M}(\beta\mathbb{S})$, where
$\beta\mathbb{S}$ denotes the Stone-{\v C}ech compactification of
$\mathbb{S}$.  A direct analog of our observability condition in this
setting would thus require that no two initial measures give rise to the
same observation statistics, even when those measures have some mass
distribution ``at infinity''. Though it is perhaps not surprising that
the observability ``at infinity'' plays a role in this setting, the
space $\beta\mathbb{S}$ is sufficiently unwieldy that a direct extension
of this type does not appear to lead to a useful theory.

In the remainder of this section we discuss two simple extensions of our
results to the non-compact case.  The first approach is inspired by the
previous discussion; if the signal admits a \textit{tractable}
compactification $\alpha\mathbb{S}$, our previous results can be
applied.  This yields stability of those functions in
$\mathcal{C}_b(\mathbb{S})$ which admit a continuous extension to
$\alpha\mathbb{S}$.  The second approach assumes that the signal process
is tight, so that the difficulties of a transient signal are avoided. In
this case it is no longer necessary that $\mathcal{O}$ is the uniform
closure of $\mathcal{O}^0$; using tightness, it is sufficient to
consider the closure of $\mathcal{O}^0$ with respect to the topology of
uniform convergence on compact sets.  This resolves our problems, as the
dual of $\mathcal{C}_b(\mathbb{S})$ endowed with the latter topology is
$\mathcal{M}_c(\mathbb{S})$, the space of compactly supported finite
signed measures.

\begin{rem}
It is only fair to remark that neither of these approaches is
particularly satisfying.  In particular, the natural test case for the
theory, the Kalman filter with an unstable signal, is not covered by
these approaches. (The Kalman filter for stable signals is not
particularly interesting, as such filters are always stable regardless
of observability; see, e.g., the result in the appendix).  Further work
is needed to develop an approach that covers unstable signals in a more
satisfactory manner. 
\end{rem}

\subsection{Compactification}

We consider a signal-observation model $(X_t,Y_t)$ as in section
\ref{sec:model}, except that $\mathbb{S}$ is not assumed to be compact. 
Let us assume, furthermore, that the observations are of the white noise
or counting type as in section \ref{sec:obsmodels}, and that the
observation function is continuous and bounded
$h\in\mathcal{C}_b(\mathbb{S})$.

\begin{defn}
Let $\alpha\mathbb{S}$ be a compact Polish space, and consider a
filtering model $(X^\alpha_t,Y^\alpha_t)$ with signal state space
$\alpha\mathbb{S}$ and observation function $h^\alpha$ (the observation
model is chosen to coincide with that of $(X_t,Y_t)$). Then
$(X_t^\alpha,Y_t^\alpha)$ is a \textit{compactification} of $(X_t,Y_t)$
if there exists a continuous injection
$\pi:\mathbb{S}\to\alpha\mathbb{S}$ such that
\begin{enumerate}
\item $h^\alpha\in\mathcal{C}_b(\alpha\mathbb{S})$ and $h=h^\alpha\circ\pi$;
\item For any $\mu\in\mathcal{P}(\mathbb{S})$, the process 
$(\pi(X_t),\pi(Y_t))$ with initial law $(X_0,Y_0)\sim\mu$ has
the same law as $(X_t^\alpha,Y_t^\alpha)$ with initial law
$(X_0^\alpha,Y_0^\alpha)\sim \mu\circ\pi^{-1}$.
\end{enumerate}
The set $\alpha\mathbb{S}\backslash\pi(\mathbb{S})$ is called the set of
\textit{points at infinity}.
\end{defn}

Denote by
$\mathcal{C}_\alpha(\mathbb{S})=\{f\in\mathcal{C}_b(\mathbb{S}):
f=f^\alpha\circ\pi\mbox{ for some
}f^\alpha\in\mathcal{C}_b(\alpha\mathbb{S})\}$ the set of bounded
continuous functions on $\mathbb{S}$ that admit a continuous extension
to $\alpha\mathbb{S}$ (note that $\mathcal{C}_\alpha(\mathbb{S})$ will
always be strictly smaller than $\mathcal{C}_b(\mathbb{S})$, unless
$\alpha\mathbb{S}=\beta\mathbb{S}$).  Then the following results follow
immediately from our definitions.

\begin{prop}
If the filter for the model $(X_t^\alpha,Y_t^\alpha)$ is stable, then
the filter for the model $(X_t,Y_t)$ is $\alpha$stable, i.e.,
$\mu\ll\nu$ implies
$$
        |\mathbf{E}^\mu(f(X_t)|\mathscr{F}_t^Y)-
        \mathbf{E}^\nu(f(X_t)|\mathscr{F}_t^Y)|\to 0
        \quad
        \mathbf{P}^\mu\mbox{-a.s.}
        \quad
        \mbox{for all }
        f\in\mathcal{C}_\alpha(\mathbb{S}).
$$
\end{prop}

\begin{cor}
Observability of the filtering model $(X_t^\alpha,Y_t^\alpha)$ implies that
the filter for $(X_t,Y_t)$ is $\alpha$stable.  In particular, 
$\alpha$stability is guaranteed if $h^\alpha$ is one-to-one.
\end{cor}

We develop further a particular setting in which this result can be
exploited.  Set $\mathbb{S}=\mathbb{R}^n$, and consider a signal
which solves the It\^o stochastic differential equation 
$$
	dX_t = f(X_t)\,dt + g(X_t)\,dW_t. 
$$ 
We assume that $f,g$ are continuously differentiable and
of sublinear growth:
$$
	f,g\in C^1,\qquad
	\|f(x)\|\le K(1+\|x\|^\alpha),\quad
	\|g(x)\|\le K(1+\|x\|^\alpha),\quad
	\alpha<1.
$$
We now consider a compactification which adjoins to $\mathbb{R}^n$ a
sphere at infinity $S^{n-1}$; this can be done, for example, by choosing
$\alpha\mathbb{S}$ to be the closed unit ball $\{x\in\mathbb{R}^n:
\|x\|\le 1\}$ and setting $\pi(x)=(1+\|x\|^2)^{-1/2}x$.  E.g., when $n=1$, this
reduces to the two-point compactification $[-\infty,\infty]$ of the real
line $\mathbb{R}=\mbox{}]{-\infty},\infty[\mbox{}$, in which case
$\mathcal{C}_\alpha(\mathbb{R})$ is precisely the set of functions 
$f\in\mathcal{C}_b(\mathbb{R})$ such that $\lim_{x\to\pm\infty}f(x)$ exist.

Now choose $h\in\mathcal{C}_\alpha(\mathbb{S})$, and let $Y_t$ be a
white noise type or counting observation model with observation function
$h$.  Then it follows from \cite[example 3]{Li94} that there is a
compactification $(X_t^\alpha,Y_t^\alpha)$ of $(X_t,Y_t)$ with the
additional property that if $X_0^\alpha\in\alpha\mathbb{S}\backslash
\pi(\mathbb{S})$ a.s., then $X_t^\alpha=X_0^\alpha$ for all $t\ge 0$
a.s.\ (i.e., the points at infinity are fixed points for the
compactified signal $X_t^\alpha$).  We can exploit the latter to give a
criterion for $\alpha$stability in terms of properties of the
non-compactified model.

\begin{prop}
Suppose that the following conditions hold:
\begin{enumerate}
\item $(X_t,Y_t)$ is observable (no two initial measures give rise to the
same law of the observation process);
\item The restriction of $h^\alpha$ to $\alpha\mathbb{S}\backslash
\pi(\mathbb{S})$ is one-to-one;
\item $X_t$ does not possess an invariant manifold that is contained
in some level set $\{x\in\mathbb{S}:h(x)=u\}$ with $u\in \{h^\alpha(x):
x\in\alpha\mathbb{S}\backslash\pi(\mathbb{S})\}$.
\end{enumerate}
Then the filter for $(X_t,Y_t)$ is $\alpha$stable.  
\end{prop}

\begin{rem}
For condition (3) to be satisfied, it is sufficient to establish that we 
have
$\{h^\alpha(x):x\in\alpha\mathbb{S}\backslash\pi(\mathbb{S})\}
\cap\{h(x):x\in\mathbb{S}\}=\varnothing$.  A sufficient condition for
(1)--(3) to be satisfied is that $h^\alpha$ is one-to-one.
\end{rem}

\begin{proof}
Let $\mu,\nu$ be two measures on $\alpha\mathbb{S}$ with $\mu\ne\nu$. 
To prove the proposition, it suffices to establish that the law of
$h^\alpha(X_t^\alpha)$ is different for $X_0^\alpha\sim\mu$ and
$X_0^\alpha\sim\nu$.

First, suppose that $\mu,\nu$ are supported on
$\alpha\mathbb{S}\backslash \pi(\mathbb{S})$. Recall that if
$X_0^\alpha\in\alpha\mathbb{S}\backslash\pi(\mathbb{S})$, then
$X_t^\alpha=X_0^\alpha$ for all $t\ge 0$.  As $h^\alpha$ is one-to-one
on $\alpha\mathbb{S}\backslash\pi(\mathbb{S})$, this implies the claim.

Next, suppose that $\mu$ is supported on $\pi(\mathbb{S})$, while $\nu$
is supported on $\alpha\mathbb{S}\backslash \pi(\mathbb{S})$. If
$X_0^\alpha\sim\mu$ and $X_0^\alpha\sim\nu$ were to give rise to the
same law for $h^\alpha(X_t^\alpha)$, then $h^\alpha(X_t^\alpha)$ would
have to be an a.s.\ constant process under both measures, so that in
particular $\mu$ must be supported on the union of the invariant
manifolds of $X_t$ which are contained in level sets of $h$.  But by the
third assumption the process $h^\alpha(X_t^\alpha)$ can then not take
the same values under $\mu$ and $\nu$, so that we have a contradiction.

Now let $\mu$ and $\nu$ be arbitrary. Then we can write
$\mu=a\mu_1+(1-a)\mu_2$ and $\nu=b\nu_1+(1-b)\nu_2$, where
$a,b\in[0,1]$, $\mu_1,\nu_1$ are supported on $\pi(\mathbb{S})$, and
$\mu_2,\nu_2$ are supported on $\alpha\mathbb{S}\backslash
\pi(\mathbb{S})$.   Note that by our assumptions, the probability that
$h^\alpha(X_t^\alpha)$ is a constant process which takes values in
$\{h^\alpha(x):x\in\alpha\mathbb{S}\backslash\pi(\mathbb{S})\}$ is
precisely $1-a$ when $X_0^\alpha\sim\mu$ and $1-b$ when
$X_0^\alpha\sim\nu$.  Hence if $\mu,\nu$ give rise to the same
observation law, then $a=b$ and $\mu_2=\nu_2$ (the latter follows as
$h^\alpha$ is assumed to be one-to-one on $\alpha\mathbb{S}\backslash
\pi(\mathbb{S})$).  But then we can conclude that $\mu,\nu$ give rise 
to the same observation law only if $\mu_1$ and $\nu_1$ give rise to the
same observation law, and the latter implies $\mu_1=\nu_1$ by our first
assumption.  Hence the proof is complete.
\end{proof}

\begin{rem}
It is likely that the compactification approach described in this section
can be generalized to a larger class of signals.  However, the requirement
that $h\in\mathcal{C}_\alpha(\mathbb{S})$ appears to rule out any model
in which the observation function is unbounded.  In practice, on the other
hand, the unbounded observation case is much more natural when $\mathbb{S}$ 
is non-compact, and one would even expect stability to improve in this 
setting.  The restriction to bounded observation functions is thus a 
significant drawback of the compactification approach. In particular, this
rules out the application of this approach to the Kalman filter.
\end{rem}

\subsection{Tight signal}

In this section, we consider a signal-observation model $(X_t,Y_t)$ as
in section \ref{sec:model}, except that $\mathbb{S}$ is assumed to
be only locally compact.  Throughout this subsection, the space of
bounded continuous functions $\mathcal{C}_b(\mathbb{S})$ will always
be endowed with the topology of uniform convergence on compact sets.
This is a locally convex topology, and gives rise to the duality
$\mathcal{C}_b(\mathbb{S})^*=\mathcal{M}_c(\mathbb{S})$
\cite[proposition IV.4.1]{Con85}.  We note that lemma \ref{lem:annihilator}
extends also to this setting:
\begin{lem}
\label{lem:annihilator2}
Let $M\subset\mathcal{C}_b(\mathbb{S})$ be a linear subspace.  Then 
$(M^\perp)^\perp=\overline{M}$, where $\overline{M}$ is the closure of $M$ 
in the topology of uniform convergence on compact sets.
\end{lem}
This can be proved in the same way as lemma \ref{lem:annihilator}, or follows
as a special case of \cite[theorem V.1.8]{Con85}.  In this setting, we
will define
\begin{equation*}
\begin{split}
	\mathcal{N} &= \{\alpha\mu_1-\alpha\mu_2\in\mathcal{M}_c(\mathbb{S}):
	\alpha\in\mathbb{R},~\mu_1,\mu_2\in\mathcal{P}_c(\mathbb{S}),~
	\mu_1\eqy\mu_2\},\\
	\mathcal{O}&=\{f\in\mathcal{C}_b(\mathbb{S}):
		\mu_1(f)=\mu_2(f)\mbox{ for all }\mu_1\eqy\mu_2
		\mbox{ with }\mu_1,\mu_2\in\mathcal{P}_c(\mathbb{S})\}.
\end{split}
\end{equation*}
That is, we will consider bounded functions (which do not necessarily vanish
at infinity), but we need only consider measures which are compactly supported.
We now obtain the following analog of proposition \ref{prop:dense}.

\begin{prop}
\label{prop:dense2}
Let $\mathcal{O}^0$ be the linear span of functions of the form
$$
        \mathbf{E}_{(x,y)}(f_1(Y_{t_1}-Y_0)f_2(Y_{t_2}-Y_0)\cdots
        f_n(Y_{t_n}-Y_0)),
$$
for all $n<\infty$, $t_i\in D$ and $f_i\in\mathcal{C}_b(\mathbb{O})$,
where $D$ is a dense subset of $[0,\infty[\mbox{}$.  Then
$\mathcal{O}^0$ is dense in $\mathcal{O}$ in the topology of uniform
convergence on compact sets.
\end{prop}

The proof is identical to that of proposition \ref{prop:dense}, so we do
not repeat it.  We are now in the position to prove a stability theorem,
similar to theorem \ref{thm:stab}, provided we assume that the signal
is tight. The notion of stability is also somewhat weaker than that of
theorem \ref{thm:stab}, as a.s.\ convergence is replaced by convergence
in $L^1$.

\begin{thm}
Let $\mu\ll\nu$, and assume that the signal $X_t$ is tight in the sense
that for every $\varepsilon>0$, there is a compact set
$K_\varepsilon\subset\mathbb{S}$ such that $\mathbf{P}^\nu(X_t\in
K_\varepsilon)>1-\varepsilon$ for all $t\ge 0$.  Then
$
        \mathbf{E}^\mu(|\mathbf{E}^\mu(f(X_t)|\mathscr{F}_t^Y)-
        \mathbf{E}^\nu(f(X_t)|\mathscr{F}_t^Y)|)\xrightarrow{t\to\infty} 0
$
for any $f\in\mathcal{O}$.
\end{thm}

\begin{proof}
For $f\in\mathcal{O}^0$, the result follows directly from lemma
\ref{lem:cvgsimple}.  Now choose $\varepsilon,\delta>0$.  By the
tightness assumption and \cite[page 146(b)]{Wil91}, we can choose
a compact set $K$ such that $\mathbf{P}^\mu(X_t\in K)>1-\varepsilon$
and $\mathbf{P}^\nu(X_t\in K)>1-\varepsilon$ for all $t\ge 0$.
We also choose $f_\delta\in\mathcal{O}^0$ be such that
$\sup_{x\in K}|f(x)-f_\delta(x)|<\delta$.  Then
\begin{multline*}
        |\pi_t^\mu(f)-\pi_t^\nu(f)| \le
        |\pi_t^\mu(f)-\pi_t^\mu(fI_{K})| +
        |\pi_t^\mu(fI_{K})-\pi_t^\mu(f_\delta I_{K})| \\ +
        |\pi_t^\mu(f_\delta I_{K})-\pi_t^\mu(f_\delta)| +
        |\pi_t^\mu(f_\delta)-\pi_t^\nu(f_\delta)| +
        |\pi_t^\nu(f_\delta)-\pi_t^\nu(f_\delta I_{K})| \\
        +
        |\pi_t^\nu(f_\delta I_{K})-
                \pi_t^\nu(fI_{K})| +
        |\pi_t^\nu(fI_{K})-\pi_t^\nu(f)|,
\end{multline*}
where we have written 
$\pi_t^\mu(f)=\mathbf{E}^\mu(f(X_t)|\mathscr{F}_t^Y)$.
Note that
\begin{equation*}
\begin{split}
        &\mathbf{E}^\mu(|\pi_t^\mu(f)-\pi_t^\mu(fI_{K})|) \le
        \|f\|\,\mathbf{P}^\mu(X_t\in K^c)\le
        \varepsilon\|f\|,\\
        &\mathbf{E}^\mu(|\pi_t^\mu(f_\delta I_{K})-
        \pi_t^\mu(f_\delta)|) \le
        \|f_\delta\|\,\mathbf{P}^\mu(X_t\in K^c)\le\varepsilon
        \|f_\delta\|,
\end{split}
\end{equation*}
while
$$
        \mathbf{E}^\mu(|\pi_t^\mu(fI_{K})-\pi_t^\mu(f_\delta
                I_{K})|) \le
        \delta,\qquad
        \mathbf{E}^\mu(|\pi_t^\nu(f_\delta I_{K})-
                \pi_t^\nu(fI_{K})|)\le \delta.
$$
Now fix $\gamma>0$, and note that
\begin{multline*}
        \mathbf{E}^\mu(|\pi_t^\nu(fI_{K})-\pi_t^\nu(f)|) \le
        \|f\|\,\mathbf{E}^\nu\!\left(\frac{d\mu}{d\nu}(X_0)\,
                \pi_t^\nu(I_{K^c})\right) \\
        \le
        \|f\|\,\mathbf{E}^\nu\!\left(\frac{d\mu}{d\nu}(X_0)\,
                I_{d\mu/d\nu(X_0)>\gamma}\right)
        +
        \gamma\|f\|\,\mathbf{P}^\nu(X_t\in K^c) \\
        \le
        \|f\|\,\mathbf{E}^\nu\!\left(\frac{d\mu}{d\nu}(X_0)\,
                I_{d\mu/d\nu(X_0)>\gamma}\right)
        +
        \varepsilon\gamma\|f\|,
\end{multline*}
and similarly
$$
        \mathbf{E}^\mu(|\pi_t^\nu(f_\delta)-
                \pi_t^\nu(f_\delta I_{K})|)
        \le
        \|f_\delta\|\,\mathbf{E}^\nu\!\left(\frac{d\mu}{d\nu}(X_0)\,
                I_{d\mu/d\nu(X_0)>\gamma}\right)
        + \varepsilon\gamma\|f_\delta\|.
$$
It follows that
\begin{multline*}
        \limsup_{t\to\infty}
        \mathbf{E}^\mu(|\pi_t^\mu(f)-\pi_t^\nu(f)|)  \\ \le
        2\delta + \left(\varepsilon+\varepsilon\gamma
                + \mathbf{E}^\nu\left(\frac{d\mu}{d\nu}(X_0)\,
                        I_{d\mu/d\nu(X_0)>\gamma}\right)
        \right)(\|f\|+\|f_\delta\|).
\end{multline*}
But $\delta,\varepsilon,\gamma>0$ were arbitrary, so the result follows
by letting $\delta,\varepsilon\to 0$, $\gamma\to\infty$.
\end{proof}

\begin{rem}
If we consider all continuous functions $\mathcal{C}(\mathbb{S})$, rather than
bounded continuous functions $\mathcal{C}_b(\mathbb{S})$, then it is still
true that $\mathcal{C}(\mathbb{S})^*=\mathcal{M}_c(\mathbb{S})$ when
$\mathcal{C}(\mathbb{S})$ is topologized by uniform convergence on compact
sets.  Thus, in fact, we may consider the set of continuous observable 
functions
$$
	\mathcal{O}=\{f\in\mathcal{C}(\mathbb{S}):
		\mu_1(f)=\mu_2(f)\mbox{ for all }\mu_1\eqy\mu_2
		\mbox{ with }\mu_1,\mu_2\in\mathcal{P}_c(\mathbb{S})\},
$$
and it is still the case that $\mathcal{O}^0$ (which contains only
bounded continuous functions) is dense in $\mathcal{O}$.  This may be
exploited, under suitable additional restrictions on the signal process,
to prove stability of unbounded observable functions that satisfy an
appropriate growth condition.  For example, if $\mathbb{S}=\mathbb{R}^n$,
$\|X_t\|^p$ is uniformly integrable and $\|d\mu/d\nu\|_q<\infty$, then
the proof of the previous theorem is easily modified to show that observable
functions of polynomial growth of degree $k$ (depending on $p,q$) are stable.
See \cite[proposition 3.3]{CL} for a similar argument.
\end{rem}

\appendix

\section{On stability of the Kalman filter}

We have shown, for a reasonably general class of nonlinear filters, that
observability implies filter stability provided that we choose
absolutely continuous initial measures $\mu\ll\nu$.  A similar result
for linear filtering models is well known, and the stability of the
Kalman filter has been studied already for several decades (see, e.g.,
\cite{BJ68}; a more recent account can be found in \cite{OP96}).  These
results, however, do not tend to assume that $\mu\ll\nu$, while on the
other hand stabilizability is typically required in addition to
detectability.  The goal of this appendix is to illustrate that also in
the linear setting, the stabilizability condition can be disposed of if
we are willing to impose an absolute continuity requirement on the
initial conditions.  This highlights the separate roles of
stabilizability and detectability in this setting.

\begin{rem}
We will make no attempt at generality, and prove only the simplest
possible result, for the purpose of illustration, by applying readily 
available results from the literature.  Despite that the conclusion is hardly
surprising and that the proof is straightforward, the author could not find 
any such result in the literature.
\end{rem}

We consider the following linear signal-observation model:
\begin{equation*}
\begin{split}
        dX_t &= AX_t\,dt + B\,dW_t,\\
        dY_t &= CX_t\,dt + dB_t,
\end{split}
\end{equation*}
where $A$, $B$, and $C$ are $n\times n$, $n\times m$, and $p\times n$ 
matrices, respectively, and $X_0$ is Gaussian with mean $\hat X_0$ and 
covariance matrix $P_0$.  As is well known, the filtered estimate $\hat 
X_t=\mathbf{E}(X_t|\mathscr{F}_t^Y)$ and covariance 
$P_t=\mathbf{E}((X_t-\hat X_t)(X_t-\hat X_t)^*)$ satisfy
\begin{equation*}
\begin{split}
        d\hat X_t &= A\hat X_t\,dt + P_tC^*(dY_t-C\hat X_t\,dt),
        \\
        \frac{dP_t}{dt} &= AP_t+P_tA^* + BB^*-P_tC^*CP_t.
\end{split}
\end{equation*}
The first equation is the Kalman filtering equation, while the second is 
the Riccati equation.  We would like to compare the solution of these 
equations with the solutions $\hat X_t'$, $P_t'$ of the same equations 
with incorrect initial conditions $\hat X_0'$, $P_0'$.

\begin{prop}
Suppose $(A,C)$ is detectable and $P_0>0$, $P_0'>0$.  Then there is a 
$P_\infty$ such that $P_t,P_t'\to P_\infty$ as $t\to\infty$, and 
$\mathbf{E}\|\hat X_t-\hat X_t'\|\to 0$ as $t\to\infty$.
\end{prop}

\begin{rem} 
It is known that under the conditions of the proposition the solution of
the Riccati equation converges to a unique limit $P_\infty$ \cite{DG92}. 
There is a key difference with the stabilizable case, however: in the
current setting, the matrix $A-P_\infty C^*C$ could be singular, while
in the stabilizable case the matrix $A-P_\infty C^*C$ is guaranteed to
be strictly negative.  The current situation is thus more subtle, and
the proof of \cite{OP96} does not immediately extend to this setting.
\end{rem}

\begin{proof}
The existence of the (nonnegative definite) matrix $P_\infty$ and the 
convergence of $P_t,P_t'$ is established in \cite{DG92} (see also 
\cite{PK97} for more recent results).  It remains to establish the second 
part of the proposition.  To this end, recall that the innovation $d\bar 
B_t=dY_t-C\hat X_t\,dt$ is a Wiener process.  We begin by writing
$$
        d(\hat X_t-\hat X_t') = (A- P_t'C^*C)(\hat X_t-\hat X_t')\,dt
                + (P_t- P_t')C^*\,d\bar B_t.
$$
Now recall that the fact that $(A,C)$ is detectable implies that there 
exists a matrix $K$ such that $A-KC$ only has eigenvalues with strictly 
negative real parts.  Fix any such $K$, define $F=KC-A$, and note that we 
can write
\begin{multline*}
        \hat X_T-\hat X_T' = e^{-FT}(\hat X_0-\hat X_0') + \mbox{} 
        \\
                \int_0^T e^{-F(T-t)}
                (K-P_t'C^*)C(\hat X_t-\hat X_t')\,dt
                + \int_0^T e^{-F(T-t)}(P_t- P_t')C^*\,d\bar B_t.
\end{multline*}
We claim that each of these terms converges to zero in $L^1$.  Let us 
consider each term individually.  The first term clearly converges to 
zero, as the eigenvalues of $F$ have strictly positive real parts.  For 
the second term, note that there are constants $c_1,\lambda>0$ such that 
$\|e^{-F(T-t)}\|\le c_1\,e^{-\lambda(T-t)}$, and that as $P_t'\to 
P_\infty$ it must be the case that $\|P_t'\|$ is bounded from above by 
some constant $c_2$.  We can thus estimate 
$$
        \left\|
                \int_0^T e^{-F(T-t)}
                (K-P_t'C^*)C(\hat X_t-\hat X_t')\,dt
        \right\| \le
                c_3\int_0^T e^{-\lambda(T-t)}
                \|C(\hat X_t-\hat X_t')\|\,dt,
$$
where we have lumped the various constants into $c_3>0$.  But by
\cite[theorem 3.1]{COC99}
$$
        \mathbf{E}\left[\int_0^\infty
                \|C(\hat X_t-\hat X_t')\|^2\,dt\right]<\infty
        \quad\Longrightarrow\quad
        \int_0^\infty
                (\mathbf{E}\|C(\hat X_t-\hat X_t')\|)^2\,dt<\infty.
$$
Hence the second term converges to zero in $L^1$ as $T\to\infty$ by lemma 
\ref{lem:cvgexpw} below.  It remains to deal with the third term.  To this 
end, note that
$$
        \mathbf{E}\left(
        \left\|\int_0^T e^{-F(T-t)}(P_t- P_t')C^*\,d\bar B_t\right\|^2
        \right)
        =
        \int_0^T {|||}e^{-F(T-t)}(P_t- P_t')C^*{|||}^2\,dt,
$$
where $|||\cdot|||$ is the Frobenius norm.  It is easily seen, again using 
lemma \ref{lem:cvgexpw} below, that this expression converges to zero as 
$T\to\infty$.  Hence the stochastic integral converges to zero in $L^2$, 
and thus also in $L^1$.  The proof is complete.
\end{proof}

The following simple lemma was used in the proof.

\begin{lem}
\label{lem:cvgexpw}
Let $\lambda>0$ and $f:[0,\infty[\mbox{}\to\mathbb{R}$.  Then
$$
        \int_0^\infty f(t)^2\,dt < \infty\qquad
        \Longrightarrow\qquad
        \int_0^T e^{-\lambda (T-t)}f(t)\,dt\xrightarrow{T\to\infty}0.
$$
The result also holds if we require $f(t)\to 0$ instead of the $L^2$ 
bound.
\end{lem}

\begin{proof}
Note that for any $s>0$, we can estimate
\begin{multline*}
        \limsup_{T\to\infty}
        \left|\int_0^T e^{-\lambda (T-t)}f(t)\,dt\right| \\ \le
        \limsup_{T\to\infty}
        \left|\int_0^s e^{-\lambda (T-t)}f(t)\,dt\right| +
        \limsup_{T\to\infty}
        \left|\int_s^T e^{-\lambda (T-t)}f(t)\,dt\right|.
\end{multline*}
Clearly the first term on the right is zero.  For the second term, note 
that
$$
        \left|\int_s^T e^{-\lambda (T-t)}f(t)\,dt\right| \le
        \left[
        \int_s^T e^{-2\lambda (T-t)}\,dt
        \int_s^T f(t)^2\,dt
        \right]^{1/2},
$$
by Cauchy-Schwarz.  The result follows by letting $T\to\infty$, then 
$s\to\infty$.  To prove that the result also holds if $f(t)\to 0$, it 
suffices to repeat the proof with
$$
        \left|\int_s^T e^{-\lambda (T-t)}f(t)\,dt\right| \le
        \left[\sup_{t\ge s}|f(t)|\right]\int_s^T e^{-\lambda (T-t)}\,dt,
$$
where it should be noted that $f(t)\to 0$ implies
$\limsup_{t\to\infty}|f(t)|=0$.
\end{proof}

\vskip.2cm
\textbf{Acknowledgement.}
The author is indebted to Pavel Chigansky for countless useful discussions 
on the topic of this article.

\bibliographystyle{acm}
\bibliography{ref}

\begin{thebibliography}{10}

\bibitem{AZ97}
{\sc Atar, R., and Zeitouni, O.}
\newblock Exponential stability for nonlinear filtering.
\newblock {\em Ann. Inst. H. Poincar\'e Probab. Statist. 33\/} (1997),
  697--725.

\bibitem{Bar55}
{\sc Bartlett, M.~S.}
\newblock {\em An introduction to stochastic processes, with special reference
  to methods and applications}.
\newblock Cambridge University Press, 1955.

\bibitem{geom}
{\sc Basile, G., and Marro, G.}
\newblock {\em Controlled and Conditioned Invariants in Linear Systems Theory}.
\newblock Prentice Hall, 1992.

\bibitem{BCL}
{\sc Baxendale, P., Chigansky, P., and Liptser, R.}
\newblock Asymptotic stability of the {W}onham filter: {E}rgodic and nonergodic
  signals.
\newblock {\em SIAM J. Control Optim. 43\/} (2004), 643--669.

\bibitem{Bil}
{\sc Billingsley, P.}
\newblock {\em Convergence of Probability Measures}, second~ed.
\newblock Wiley, 1999.

\bibitem{BJ68}
{\sc Bucy, R.~S., and Joseph, P.~D.}
\newblock {\em Filtering for stochastic processes with applications to
  guidance}.
\newblock Wiley, 1968.

\bibitem{ChSurv}
{\sc Chigansky, P.}
\newblock Stability of nonlinear filters: A survey, 2006.
\newblock Mini-course lecture notes, Petropolis, Brazil. Available at {\tt
  http://www.wisdom.weizmann.ac.il/$\sim$pavel/}.

\bibitem{CL}
{\sc Chigansky, P., and Liptser, R.}
\newblock On a role of predictor in the filtering stability.
\newblock {\em Electr. Commun. Probab. 11\/} (2006), 129--140.

\bibitem{COC99}
{\sc Clark, J. M.~C., Ocone, D.~L., and Coumarbatch, C.}
\newblock Relative entropy and error bounds for filtering of {M}arkov
  processes.
\newblock {\em Math. Control Signals Systems 12\/} (1999), 346--360.

\bibitem{Con85}
{\sc Conway, J.~B.}
\newblock {\em A course in functional analysis}.
\newblock Springer-Verlag, 1985.

\bibitem{DG92}
{\sc {De Nicolao}, G., and Gevers, M.}
\newblock Difference and differential {R}iccati equations: a note on the
  convergence to the strong solution.
\newblock {\em IEEE Trans. Autmat. Control 37\/} (1992), 1055--1057.

\bibitem{DM}
{\sc Dellacherie, C., and Meyer, P.-A.}
\newblock {\em Probabilities and Potential B}.
\newblock North-Holland, 1982.

\bibitem{delyon}
{\sc Delyon, B., and Zeitouni, O.}
\newblock Lyapunov exponents for filtering problems.
\newblock In {\em Applied Stochastic Analysis\/} (1991), Gordon and Breach,
  pp.~511--521.

\bibitem{EK86}
{\sc Ethier, S.~N., and Kurtz, T.~G.}
\newblock {\em Markov processes: Characterization and convergence}.
\newblock Wiley, 1986.

\bibitem{obs3}
{\sc Gurvits, L., and Ledoux, J.}
\newblock Markov property for a function of a markov chain: A linear algebra
  approach.
\newblock {\em Lin. Alg. Appl. 404\/} (2005), 85--117.

\bibitem{Hasm60}
{\sc {Has'minski\u\i}, R.~Z.}
\newblock Ergodic properties of recurrent diffusion processes and stabilization
  of the solution of the {C}auchy problem for parabolic equations.
\newblock {\em Teor. Verojatnost. i Primenen. 5\/} (1960), 196--214.

\bibitem{obs1}
{\sc Ito, H., Amari, S.-I., and Kobayashi, K.}
\newblock Identifiability of hidden markov information sources and their
  minimum degrees of freedom.
\newblock {\em IEEE Trans. Inf. Th. 38\/} (1992), 324--333.

\bibitem{kaijser}
{\sc Kaijser, T.}
\newblock A limit theorem for partially observed markov chains.
\newblock {\em Ann. Probab. 3\/} (1975), 677--696.

\bibitem{Kal97}
{\sc Kallenberg, O.}
\newblock {\em Foundations of modern probability}.
\newblock Probability and its Applications (New York). Springer-Verlag, 1997.

\bibitem{reeds}
{\sc Kochman, F., and Reeds, J.}
\newblock A simple proof of {K}aijser's unique ergodicity result for hidden
  markov $\alpha$-chains.
\newblock {\em Ann. Appl. Probab. 16\/} (2006), 1805--1815.

\bibitem{Kun71}
{\sc Kunita, H.}
\newblock Asymptotic behavior of the nonlinear filtering errors of {M}arkov
  processes.
\newblock {\em J. Multivar. Anal. 1\/} (1971), 365--393.

\bibitem{Kun78}
{\sc Kunita, H.}
\newblock Supports of diffusion processes and controllability problems.
\newblock In {\em Proceedings of the International Symposium on Stochastic
  Differential Equations (Res. Inst. Math. Sci., Kyoto Univ., Kyoto, 1976)\/}
  (1978), Wiley, pp.~163--185.

\bibitem{obs2}
{\sc Larget, B.}
\newblock A canonical representation for aggregated markov processes.
\newblock {\em J. Appl. Probab. 35\/} (1998), 313--324.

\bibitem{Li94}
{\sc Li, X.-M.}
\newblock Properties at infinity of diffusion semigroups and stochastic flows
  via weak uniform covers.
\newblock {\em Potential Anal. 3\/} (1994), 339--357.

\bibitem{LS01a}
{\sc Liptser, R.~S., and Shiryaev, A.~N.}
\newblock {\em Statistics of Random Processes {I}. {G}eneral Theory},
  second~ed.
\newblock Springer, 2001.

\bibitem{OP96}
{\sc Ocone, D., and Pardoux, E.}
\newblock Asymptotic stability of the optimal filter with respect to its
  initial condition.
\newblock {\em SIAM J. Control Optim. 34\/} (1996), 226--243.

\bibitem{PK97}
{\sc Park, P., and Kailath, T.}
\newblock Convergence of the {DRE} solution to the {ARE} strong solution.
\newblock {\em IEEE Trans. Autmat. Control 42\/} (1997), 573--578.

\bibitem{rao}
{\sc Rao, M.}
\newblock On modification theorems.
\newblock {\em Trans.\ Amer. Math. Soc. 167\/} (1972), 443--450.

\bibitem{RY98}
{\sc Revuz, D., and Yor, M.}
\newblock {\em Continuous Martingales and {B}rownian Motion}, third~ed.
\newblock Springer, 1999.

\bibitem{Rud73}
{\sc Rudin, W.}
\newblock {\em Functional Analysis}.
\newblock McGraw-Hill, 1973.

\bibitem{Trig76}
{\sc Triggiani, R.}
\newblock Extensions of rank conditions for controllability and observability
  to {B}anach spaces and unbounded operators.
\newblock {\em SIAM J. Control Optim. 14\/} (1976), 313--338.

\bibitem{Wil91}
{\sc Williams, D.}
\newblock {\em Probability with martingales}.
\newblock Cambridge University Press, 1991.

\end{thebibliography}

\end{document}